\documentstyle[12pt]{article}
\input {cyracc.def}
\def\baselinestretch{1.1} \topmargin -12pt
\font\tencyr=wncyr7 scaled \magstep1
\def\rus{\tencyr\cyracc}

\newcommand{\re}[1]{\mbox{\rm  (\ref{#1})}}

\newenvironment{proof}
{\noindent {\sl Proof.}\quad }{\hfill 
$\square$ \vskip1.1ex\noindent }

\newenvironment{proof*}
{\noindent {\sl Proof.}\quad }{\hfill 
$\square$}

\renewcommand{\theequation}{\thesection .\arabic{equation}}

\renewcommand{\thesubsubsection}{\theequation .\arabic{subsubsection}}

\catcode`\@=\active
\catcode`\@=11
\def\@eqnnum{\hbox to .01pt{}\rlap{\bf \hskip -\displaywidth(\theequation)}}
\catcode`\@=12

\newenvironment{s}[1]
{ \vskip1.2ex \refstepcounter{equation}
\noindent {\bf \theequation\quad #1.} \begin{sl}}{\end{sl}
\vskip1.1ex\noindent }

\newenvironment{rem}[1]
{ \vskip1.2ex \refstepcounter{equation}
\noindent {\bf \theequation\quad {#1}.} }{ \vskip1.1ex\noindent }

\newenvironment{subs}[1]
{\vskip1.2ex \refstepcounter{equation}
\noindent {\bf (\theequation)\quad #1.} }{\quad}

\newcommand {\sekt}[1]
{{\vskip2.5ex\refstepcounter{section}\setcounter{equation}{0}
\noindent\large \bf \thesection\quad \parbox[t]{424pt}{#1}
\nopagebreak\vskip1.5ex\noindent}}

\newcommand {\ah}{{\frak a}}

\newcommand {\g}{{\frak g}}
\newcommand {\h}{{\frak h}}

\newcommand {\me}{{\frak m}}

\newcommand {\te}{{\frak t}}
\newcommand {\tri}{{\frak sl}_2}
\newcommand {\ut}{{\frak u}}

\newcommand {\esi}{\varepsilon}
\newcommand {\ap}{\alpha}

\newcommand {\Lb}{\Lambda}
\newcommand {\lb}{\lambda}
\newcommand {\vp}{\varphi}
\newcommand {\OT}{\ov{\Theta}}

\newcommand {\V}{{\Bbb V}}
\newcommand {\W}{{\Bbb W}}

\newcommand {\cC}{{\cal C}}

\newcommand {\ad}{{\mathrm{ad\,}}}

\newcommand {\ch}{{\mathrm{ch\,}}}

\newcommand {\dg}{{\mathrm{deg\,}}}

\newcommand {\Ima}{{\mathrm{Im\,}}}
\newcommand {\rk}{{\mathrm{rk\,}}}

\newcommand {\GR}[2]{{\mathrm{{\bf #1}}}_{#2}}

\newcommand {\ov}{\overline}
\newcommand {\un}{\underline}

\newcommand {\vno}[1]{\vskip#1 ex\noindent}
\newcommand {\rar}{\rightarrow}

\newcommand {\qu}{\hfill $\square$}
\newcommand {\beq}{\begin{equation}}
\newcommand {\eeq}{\end{equation}}

\newcommand{\odin}{\rm{{I}\!\!\!1}}
\newfam\Bbbfam\newfam\eufam\newfam\Eufam%
\font\olala=msam10 scaled 1200%
\font\Bbbfont=msbm10 scaled 1200%
\font\Bbbsmallfont=msbm7 scaled 1000%
\font\papa=msbm8 scaled 1000%
\textfont\Bbbfam=\Bbbfont%
\scriptfont\Bbbfam=\Bbbsmallfont
\scriptscriptfont\Bbbfam=\papa%
\font\euzw=eufm10 scaled 1200%
\font\euac=eufm7 scaled 1200%
\font\euacc=eufm5 scaled 1200%
\textfont\eufam=\euzw\scriptfont\eufam=\euac\scriptscriptfont\eufam=\euacc%
\font\Euzw=eufm10 scaled \magstep2%
\font\Euac=eufm7 scaled \magstep2%
\textfont\Eufam=\Euzw\scriptfont\Eufam=\Euac%
\def\frak{\fam\eufam}%
\def\Frak{\fam\Eufam}%
\def\Bbb{\fam\Bbbfam}%
\def\bbb{\fam\Bbbfam\scriptscriptstyle}
\def\varnothing{\hbox {\Bbbfont\char'077}}
\def\square{\hbox {\olala\char"03}}
\def\Bbbk{\hbox {\Bbbfont\char'174}}

\begin{document}
\setlength{\parskip}{2pt plus 4pt minus 0pt}

\noindent
{\Large \bf  \parbox[t]{424pt}{The exterior algebra and `Spin'
of an orthogonal $\Frak g$-module}}
\bigskip \\
{\bf Dmitri I. Panyushev\footnote{supported in part by 
R.F.F.I. Grant {\rus N0}\,98--01--00598} }

\medskip
\smallskip

\noindent{\large \bf Introduction}
\vno{2}%
Let $\g$ be a reductive algebraic Lie algebra over an algebraically closed 
field $\Bbbk$ of characteristic zero and $G$ is the corresponding connected 
and simply connected group. 
\\[.6ex]
The symmetric algebra 
of a (finite-dimensional) $\g$-module $\V$
is the algebra of polynomial functions on the dual space
$\V^*$. Therefore one can study the algebra of symmetric invariants 
using geometry of $G$-orbits in $V^*$. 
In case of the exterior algebra, $\wedge^\bullet\V$,
lack of such geometric picture results by now in absence of general structure 
theorems for the algebra of skew-invariants $(\wedge^\bullet\V)^\g$.
One may find in the literature several interesting results related to 
skew-symmetric invariants. We only mention
Kostant's computation for cohomology of the nilradical of 
a parabolic subalgebra in $\g$ \cite{cohom} and
R.\,Howe's classification of ``skew-multiplicity-free'' $\g$-modules
\cite[ch.\,IV]{howe}.
But the general situation still remains unsatisfactory, and developing of
Invariant Theory in the skew-symmetric setting represents an attractive problem.
\\[.6ex]
In this paper, we begin with describing all irreducible
orthogonal $\g$-modules such that
$(\wedge^\bullet\V)^\g$ is again an exterior algebra.
It is shown that in this case either $\V\simeq\g$ and hence $\g$ is simple
or $\g\oplus\V$ has a structure of {\sl simple\/} ${\Bbb Z}_2$-graded
Lie algebra, which quickly leads to a short classification, see
Table~\ref{final}.
Obviously, none of the symplectic representations
(with $\dim\V > 2$) can have an exterior algebra of skew-invariants. 
But the situation for the representations of ``general
type'' is not yet clear.
\\[.6ex]
In case $\V$ is orthogonal, a better understanding of the $\g$-module
structure of $\wedge^\bullet\V$ can be achieved through the notion of
{\it `Spin'} of $\V$. This goes as follows. Let $\pi:\g\rar {\frak so}(\V)$ 
be the corresponding representation. Restricting the spinor representation
of ${\frak so}(\V)$ to $\g$ gives us a $\g$-module, which is denoted by
$Spin(\V)$. The motivation came from Kostant's result that 
$Spin(\g)$ is a primary $\g$-module;
namely, $Spin(\g)=2^{[\mathrm{\footnotesize \rk}\g/2]}\V_\rho$, the 
highest weight
$\rho$ being the half-sum of the positive roots \cite[p.\,358]{cohom}. 
The main property of $Spin(\V)$ is that, depending on parity of
$\dim\V$, $\wedge^\bullet\V$ is isomorphic to either $Spin(\V)^{\otimes 2}$
or $2{\cdot}Spin(\V)^{\otimes 2}$. It is thus interesting to find the 
orthogonal $\g$-modules, where $Spin(\V)$ has a simple structure.
In general, $Spin(\V)$, as element of the representation ring, has a numerical
factor depending on the zero-weight multiplicity. Omitting this factor
yields a $\g$-module, which is called the {\it reduced `Spin'} of 
$\V$ and denoted
by $Spin_0(\V)$; e.g. $Spin_0(\g)=\V_\rho$. In a sense, $Spin_0(\V)$ behaves
better than $Spin(\V)$. For, regardless of parity of $\dim\V$,
we have $\wedge^\bullet\V\simeq 2^{m(0)}{\cdot}Spin_0(\V)^{\otimes 2}$,
where $m(0)$ is the zero-weight multiplicity, and 
$Spin_0(\V_1\oplus\V_2)=Spin_0(\V_1)\otimes Spin_0(\V_2)$.
\\[.6ex]
An orthogonal $\g$-module $\V$ is said to be {\it co-primary\/}, if 
$Spin_0(\V)$ is irreducible. In sections~\ref{spin1} and \ref{spin2}, 
a classification of the co-primary modules is obtained. To this end, we
give a geometric description of some highest weights of $Spin_0(\V)$.
These weights are called {\it extreme\/}. The assumption that
$Spin_0(\V)$ has a unique extreme weight imposes strong constraints on 
the weight structure of $\V$. Using this, one shows that $\g$ must be simple
whenever $\V$ is an irreducible faithful co-primary $\g$-module, and
that any reducible co-primary module is being obtained by iterating the
``direct sum'' procedure: \\
\centerline{$(\g_i,\V_i),\ i=1,2 \ \mapsto
(\g_1\oplus\g_2,\V_1\oplus\V_2)$.}
It is thus sufficient to classify the irreducible co-primary modules. 
The resulting list appears to be rather short:

1. $\g$ is simple and $\V\simeq\g$;

2. $\g$ is of type $\GR{B}{n}$ or $\GR{C}{n}$ or $\GR{F}{4}$, and
$\V$ is the little adjoint module;

3. $\g={\frak so}(\W)$ and 
$\V=\{\mbox{the Cartan component of }{\cal S}^2\W\}$;
$\dim\W=3,\,5,\,7,\dots$.
\\ 
In case 2, $\g$ has roots of two lengths and the $\g$-module whose 
highest weight
is the short dominant root is called {\it little adjoint} (l.a.). 
Actually, we give 
a unified proof for the fact that the l.\,a. module is co-primary
whenever the  ratio of root lengths is $\sqrt{2}$. Note that,
for $\GR{G}{2}$, where this ratio is $\sqrt{3}$,
the l.\,a. module is not co-primary. A true reason why this is so
is that the l.\,a. module for $\GR{G}{2}$ is not the isotropy 
representation of a symmetric space, whereas this is the case for 
{\bf B, C}, and {\bf F}. Furthermore, all representations listed above
are the isotropy representations of symmetric spaces. A curious coincidence
in this regard is the following. Let $\tilde G/G$ be a symmetric
space and $\tilde\g=\g\oplus\V$ the corresponding Lie algebra decomposition.
Then the $\g$-module $\V$ is co-primary if and only if $\g$ is non-homologous
to zero in $\tilde\g$. Another by-product of our classifications is that any
co-primary module has an exterior algebra of skew-invariants.
\\[.5ex]
Having observed that any co-primary representation is a very specific isotropy
representation, one may suggest that $Spin_0(\V)$ admits a nice description
for {\sl all\/} symmetric spaces. This is really the case, and a transparent
formulation can be given for the inner involutory automorphisms.
Let $\g=\g_0\oplus\g_1$ be a ${\Bbb Z}_2$-grading of inner type,
i.e., $\rk\g_0=\rk\g$. As the $\g_0$-module $\g_1$ has no zero weight,
$Spin(\g_1)=Spin_0(\g_1)$. Choose a common Cartan subalgebra $\te$ for
$\g$ and $\g_0$, and
consider the natural inclusion of the Weyl groups
$W_0\subset W$. Although $W_0$ is not necessary a parabolic subgroup of $W$,
each coset $wW_0$ contains a unique element of minimal length 
(see \ref{coset}).
Let $W^0\subset W$ be the set of such elements.
Then the irreducible constituents of the $\g_0$-module $Spin(\g_1)$ are 
parameterized by $W^0$. Namely, $\displaystyle
Spin(\g_1)=\bigoplus_{w\in W^0}\V_{\lb_w}$, where $\lb_w=w^{-1}\rho-\rho_0$
is the highest weight, see section~\ref{spin4}. Moreover, the weights
$\lb_w$ ($w\in W^0$) are distinct and hence $Spin(\g_1)$ is a multiplicity
free $\g_0$-module. It is worth noting that the above expression for 
$Spin(\g_1)$ is equivalent to an identity for root systems that
seem to have not been observed before. Let $\Delta$ be the root system
of $(\g,\te)$ and let $\Delta^+=\Delta^+_0\sqcup\Delta^+_1$ be the partition
of the set of positive roots corresponding to the sum
$\g=\g_0\oplus\g_1$. In this situation, one can introduce the 
``cunning'' parity
$\tau:W\rar\{1,-1\}$, determined by $\Delta_0^+$. 
If $w\in W_0$, then $\tau(w)=(-1)^{l_0(w)}$, where $l_0({\cdot})$ is
the length in $W_0$ relative to the set of positive roots $\Delta_0^+$.
To extend $\tau$ to $W$, one uses the aforementioned subset $W^0$
(see section~\ref{spin3}). 
Then the identity reads
\[
\sum_{w\in W}\tau(w)e^{w\rho}=\prod_{\ap\in \Delta_0^+}
(e^{\ap/2}-e^{-\ap/2})\prod_{\mu\in \Delta_1^+}
(e^{\mu/2}+e^{-\mu/2})\ . 
\]
For the outer involutory automorphisms, the final description of $Spin_0(\g_1)$
is almost identical to the previous one, see section~\ref{spin5}. 
However, it requires much more
preparations and its proof uses the classification of
involutory automorphisms. Our arguments suggest that there should exist
interesting connections between cohomology of symmetric spaces,
twisted affine Kac--Moody algebras, and $Spin(\g_1)$.
\\[.6ex]
The description of the highest weights of $Spin(\g_1)$ (for all involutions!)
shows that these weights are extreme. This also implies the following
claim (see sect.\,7):

Let $\Phi(\ ,\ )$ be an invariant bilinear form on $\g$ and  
$\Phi(\ ,\ )_0$ its restriction to $\g_0$. Let $c_0\in U(\g_0)$ be the 
Casimir element with respect to $\Phi(\ ,\ )_0$. Then $c_0$ acts scalarly on
$Spin(\g_1)$; the value is $(\rho,\rho)-(\rho_0,\rho_0)$, where 
$(\ ,\ )$ is the $W$-invariant bilinear form on $\te^*$ induced by
$\Phi(\ ,\ )$. 
\\
A similar result holds for the isotropy representation $\h\rar {\frak so}(\me)$
of non-symmetric space $G/H$, if $\rk\h=\rk\g$ and one considers the
submodule of $Spin(\me)$ generated by the extreme weight vectors.
\\
Recently, B.\,Kostant obtained a series of nice results for $Spin(\g)$
\cite{cliff}. Since the adjoint representation is one of the isotropy
representations of symmetric spaces,
our results for $Spin(\g_1)$ suggest that many parts of Kostant's theory can
be generalized to the setting of arbitrary symmetric spaces.
\\[.8ex]
{\bf Main notation.} 
$\g$ is a reductive Lie algebra with a fixed triangular decomposition:
$\g=\ut^+\oplus\te\oplus\ut^-$. All $\g$-modules are assumed to be
finite-dimensional.\\
$\Delta$ (resp. $\Delta^+$) is the set of roots (resp. positive roots);
$\Pi\subset\Delta^+$ is the set of simple roots; 
$\Pi=\{\ap_i\}_{i\in I}$
and $\vp_i$ is the
fundamental weight corresponding to $\ap_i$. 
For simple Lie algebras, we follow the numeration of the simple roots 
from \cite{vion} and \cite{al}.\par
${\cal P}$ -- the lattice of integral weights, 
${\cal P}_+$ -- the monoid of dominant integral weights. \par
$W=N_G(\te)/Z_G(\te)=N_G(\te)/T$ -- the Weyl group; 
for $\beta\in\Delta$, $s_\beta$ is the 
reflection in $W$.\\
${\cal P}_{\Bbb Q}={\cal P}\otimes_{\Bbb Z}{\Bbb Q}\subset\te^*$ 
and $(\ ,\ )$ is
the $W$-invariant positive-definite scalar product in ${\cal P}_{\Bbb Q}$
determined by a non-degenerate invariant bilinear form $\Phi(\ ,\ )$ on $\g$.\\
If $M\subset\cal P$ is any finite set of weights, then $|M|=\sum_{m\in M}m$;
$\rho:=\frac{1}{2}|\Delta^+|$. \\
If $\lb\in {\cal P}_+$, then $\V_\lb$ stands for the irreducible
$\g$-module with highest weight $\lb$.
\\[.7ex]
{\small 
{\bf Acknowledgements.} I would like to thank A.L.\,Onishchik for 
conversations about cohomology of compact homogeneous spaces. I am indebted
to B.\,Kostant for drawing my attention to Conlon's paper \cite{conlon}.
Part of this
work was done while I was visiting Universit\'e de Poitiers. I am grateful to 
this institution for the great hospitality I enjoyed there.
}

\sekt{Orthogonal $\Frak g$-modules with an 
exterior algebra of skew-invariants\label{grassmann}}%
Let $\V$ be a $\g$-module. Study of the algebra $({\cal S}^\bullet\V)^\g$
of symmetric (or polynomial) invariants is the subject of a rich and
well-developed theory. In contrast, little is known about
the algebra $(\wedge^{\bullet}\V)^\g$ of skew-invariants. The 
skew-symmetric theory has some parallels to the symmetric case, and many
interesting differences. We begin with two observations.
\\[1ex]
\refstepcounter{equation}{\bf (\theequation)}\qquad
Put $n=\dim\V$. Suppose $\g\subset{\frak sl}(\V)$, e.g. $\g$ is semisimple. 
Then $\wedge^n\V$ is a trivial $\g$-module. Therefore
$\dim (\wedge^\bullet\V)^\g\ge 2$, $\wedge^i\V$ and $\wedge^{n-i}\V$
are isomorphic $\g$-modules, and the Poincar\'e polynomial
of $(\wedge^\bullet\V)^\g$ is symmetric.
\\[1ex]
\refstepcounter{equation}{\bf (\theequation)}\qquad \label{ext-char}
Let  ${\cal P}(\V)$ be the set of all weights of $\V$
relative to $\te\subset\g$  and
$\V=\oplus_{\mu\in{\cal P}(\V)} \V^\mu$ the weight decomposition.
Set $m(\mu)=\dim\V^\mu$.
Recall that the character of $\V$ is the element of the group algebra
${\Bbb Z}[{\cal P}]$ defined by $\ch\V=\sum_{\mu\in{\cal P}(\V)}m(\mu)e^\mu$.
Then, $t$ being an indeterminate, we have
\[
\sum_{i=1}^n (\ch\wedge^i\V) t^i=
\prod_{\mu\in{\cal P}(\V)}(1+te^\mu)^{m(\mu)} \ .
\]
In particular,
$\ch\wedge^\bullet\V=2^{m(0)}\prod_{\mu\ne 0}(1+e^\mu)^{m(\mu)}$.
Obviously, $\prod_{\mu\ne 0}(1+e^\mu)^{m(\mu)}$
must be the character of a $\g$-module, say $\W$.
Hence $1\le \dim(\wedge^\bullet\V)^\g= 2^{m(0)}{\cdot}\dim\W^\g$ and
therefore $\dim(\wedge^\bullet\V)^\g\ge 2^{m(0)}$. 
\\[1ex]
It is natural to first
describe $\g$-modules, where the algebra of skew-invariants has a simple
structure.
\\[.6ex]
{\bf Definition.} The algebra $(\wedge^\bullet\V)^\g$ is said to be 
{\it free\/} (or an {\it exterior algebra\/}), if there exists a graded 
subspace $P\subset (\wedge^\bullet\V)^\g$ such that 
$(\wedge^\bullet\V)^\g$ is the exterior algebra over $P$. 
\\[.6ex]
Suppose $\g\subset{\frak sl}(\V)$ and 
$(\wedge^\bullet\V)^\g$ is an exterior algebra. Let
$P=\langle p_1,\dots,p_l\rangle$ with $\dg p_i=d_i$. Then $0\ne p_1\wedge\dots
\wedge p_l$ must be an element of $\wedge^n\V$. Hence
$\sum_i d_i=n$. It follows from the definition that
all the $d_i$'s must be odd whenever $l>1$. However, if
$l=1$, then $d_1=n$ is allowed to be even. In other words, all 2-dimensional
algebras of skew-invariants are proclaimed to be exterior. 
\\[.8ex]
{\bf Example.} Let $\V=\g$. Then $(\wedge^\bullet\g)^\g$ is free.
Here $l=\rk\g$ and $d_i=2m_i+1$, where $m_1,\dots,m_l$ are the
exponents of $\g$. A purely algebraic proof of this result
was given by Koszul \cite{kos}.
\\[.6ex]
From now on,  $\V$ is an orthogonal $\g$-module, 
i.e., we are given a representation $\pi:\g\rar {\frak so}(\V)$. 
In particular, $\g\subset {\frak sl}(\V)$.
\begin{s}{Lemma} \label{4=0}
Let 
$\V$ be  an irreducible orthogonal $\g$-module with $\V^\g=0$.
Suppose $(\wedge^\bullet\V)^\g$ is free.
Then either $(\wedge^4\V)^\g=0$ or $\dim\V=4$.
\end{s}\begin{proof}
Since $(\wedge^1\V)^\g=(\wedge^2\V)^\g=0$, any nonzero
element of $(\wedge^4\V)^\g$ is a generator of $(\wedge^\bullet\V)^\g$.
\end{proof}%
Let $\mu: \V\times\V^*\rar \g^*\simeq\g$ be the moment mapping associated with
the standard symplectic structure on $\V\times\V^*\simeq T^*(\V)$.
Identifying $\V$ and $\V^*$, one obtains an {\it anti-commutative\/}
bilinear
mapping $\bar\mu: \V\times\V\rar\g$. Using the $\g$-invariant symmetric
bilinear forms $\Phi(\ ,\ )$ and $(\ ,\ )_\V$,  one may explicitly
define $\bar\mu$ by
\beq \label{mubar}
\Phi(\bar\mu(v_1,v_2), g):=(v_2, g{\cdot}v_1)_\V \ ,
\eeq
where $v_1,v_2\in\V$ and $g{\cdot}v_1$ is a shorthand for
$\pi(g)v_1$. This $\bar\mu$ yields an anti-commutative multiplication,
denoted by $[\ ,\ ]\widetilde{\ }$, in $\g\oplus\V$:
\[
[(g_1,v_1), (g_2,v_2)]\widetilde{\ }:=( [g_1,g_2]+\bar\mu(v_1,v_2),
g_1{\cdot}v_2-g_2{\cdot}v_1) \ .
\]
The following assertion is stated in \cite[p.\,152]{conlon}, in the context
of compact group representations, as the ``Cartan-Kostant theorem''. 
It is an easy
part of Kostant's characterization  of the isotropy representation of
compact homogeneous spaces [loc.\,cit].
\begin{s}{Proposition} \label{jacobi}
The multiplication $[\ ,\ ]\widetilde{\ }$ satisfies the Jacobi identity
if and only if the skew-symmetric $\g$-invariant 4-form on $\V$
\[
(v_1,v_2,v_3,v_4)\stackrel{\kappa}{\mapsto} 
\Phi(\bar\mu(v_1,v_2),\bar\mu(v_3,v_4))+
\Phi(\bar\mu(v_2,v_3),\bar\mu(v_1,v_4))+
\Phi(\bar\mu(v_3,v_1),\bar\mu(v_2,v_4))  
\]
is identically equal to zero.
\end{s}\begin{proof*}
By bilinearity of $[\ ,\ ]\widetilde{\ }$, it suffices to verify the
Jacobi identity for 4 sorts of triples:
(i) $(g_1,g_2,g_3)$, \quad
(ii) $(g_1,g_2,v_1)$, \quad
(iii) $(g_1,v_1,v_2)$, \quad
(iv) $(v_1,v_2,v_3)$,\\
where $g_i\in\g$ and $v_i\in\V$.
The Jacobi identity is always satisfied for cases (i)--(iii), because,
respectively, $\g$ is a Lie algebra, $\V$ is a $\g$-module, and
$\bar\mu$ is a homomorphism of $\g$-modules.
For $v_1,v_2,v_3\in \V$, the identity means that \\[.6ex]
\centerline{
$\bar\mu(v_1,v_2){\cdot}v_3 +
\bar\mu(v_2,v_3){\cdot}v_1+\bar\mu(v_3,v_1){\cdot}v_2=0\in\V$
}
or 
\[
(\bar\mu(v_1,v_2){\cdot}v_3+
\bar\mu(v_2,v_3){\cdot}v_1+\bar\mu(v_3,v_1){\cdot}v_2, v_4)_\V=0
\]
for any $v_4\in\V$. Using Eq.~\re{mubar}, one rewrites the last equality
as the condition that the mapping $\kappa: \V^{\otimes 4}\to \Bbbk$ is zero. 
It is also easily seen that $\kappa$ is skew-symmetric and $\g$-invariant.
\end{proof*}%
\begin{s}{Corollary} \label{z2}
If $(\wedge^4\V)^\g=0$, then $\g\oplus\V$, endowed with multiplication
$[\ ,\ ]\widetilde{\ }$, is a ${\Bbb Z}_2$-graded Lie algebra.  
\end{s}%
Notice that the condition $(\wedge^4\V)^\g=0$ is not necessary for
$\g\oplus\V$ to be a ${\Bbb Z}_2$-graded
Lie algebra. We are going to list all irreducible
orthogonal $\g$-modules $\V$ such
that $(\wedge^\bullet\V)^\g$ is free.
\begin{s}{Theorem}  \label{main1}
Let $\g$ be semisimple and $\V$ a faithful orthogonal irreducible
$\g$-module. Suppose $(\wedge^\bullet\V)^\g$ is free. Then either
$\g$ is simple and $\V\simeq\g$ or $\tilde\g:=\g\oplus\V$
is a \un{simple} ${\Bbb Z}_2$-graded Lie algebra.
\end{s}\begin{proof}
1. Assume that $\dim\V\ne 4$. By Lemma~\ref{4=0} and Corollary~\ref{z2},
it follows that $[\ ,\ ]\widetilde{\ }$ makes $\tilde\g$ a
${\Bbb Z}_2$-graded Lie algebra. Let $\ah\subset\tilde\g$ be an ideal.
Then $\ah\cap\V$ and $\ah\cap\g$ are $\g$-stable spaces. 

(i) If $\ah\cap\V=\V$, then $\ah$ also contains $\bar\mu(\V,\V)=
[\V,\V]\widetilde{\ }$. Since $\V$ is faithful, $\bar\mu(\V,\V)$
meets all the simple components of $\g$. Therefore $\ah=\tilde\g$.

(ii) If $\ah\cap\g\ne 0$, then, since $\V$ is faithful,
$(\ah\cap\g){\cdot}\V\ne 0$. That is, $\ah\cap\V\ne 0$ and we are back
in part (i).

(iii) If $\ah\cap\V = 0$ and $\ah\cap\g=0$, then the $\g$-module $\ah$
is isomorphic to its projections to both $\V$ and $\g$. Hence
$pr_\g(\ah)\simeq\ah\simeq pr_\V(\ah)\simeq \V$. Therefore
$pr_\g(\ah)$ is a a simple component of $\g$. 
As $\V$ is a faithful $\g$-module, we conclude that $\g\simeq\V$ and
therefore $\g$ is simple in this case. 
Here $\tilde\g$ is the sum of two isomorphic ideals, 
$\tilde\g\simeq\ah\oplus\ah$.
The subalgebra $\g$, which is isomorphic to $\ah$, 
is the diagonal in
$\tilde\g$, and $\V=\{(x,-x)\mid x\in\ah\}$.
\\
2. Assume that $\dim\V=4$. Then $\g\subset {\frak so}_4={\frak so}(\V)$.
Obviously, ${\frak so}_4\oplus\V\simeq {\frak so}_5$, and one easily verifies
that $(\wedge^\bullet\V)^\g$ is not free for any proper reductive
subalgebra $\g$ of ${\frak so}_4$.
\end{proof}%
As is mentioned above, $(\wedge^\bullet\g)^\g$ is free. Thus, all other
irreducible orthogonal modules with free algebra of skew-invariants
arise in connection with 
${\Bbb Z}_2$-gradings of {simple}
Lie algebras. If $\tilde\g=\g\oplus\V$ is a simple ${\Bbb Z}_2$-graded
Lie algebra, then $(\wedge^\bullet\V)^\g$ is isomorphic to $H^*(\tilde G/G)$,
the cohomology ring of the symmetric space 
$\tilde G/G$ \cite[\S\,9,\,n.11]{al}. 
The cases, where $H^*(\tilde G/G)$ is an exterior
algebra, are well known, see \cite[\S\,13,\,Th.\,1]{al}. Note however that
our interpretation of ``exterior algebras'' is a bit wider. In case
$\dim H^*(\tilde G/G)=2$, the generator is allowed to be of even
degree. The resulting classification is presented in Table~\ref{final}.

\begin{table}[htbp]
\caption{The irreducible orthogonal representations with free algebra 
of skew-invariants} \label{final}
\vskip1ex\centerline
{
\begin{tabular}{clcc}
$\g$ & $\V$  & $\dim P$ &  Poincar\'e polynomial of $(\wedge^\bullet\V)^\g$ 
\\ \hline 
{\renewcommand{\baselinestretch}{.8}
\parbox{60pt}{\small any simple\\ Lie algebra}}\rule{0pt}{4ex}
& $\g$ & $\rk\g$ & $\prod_{i=1}^{{\mathrm{\scriptsize rk\,}}\g}(1+t^{2m_i+1})$
\\[2ex] 
${\frak sp}_{2n}$ ($n\ge 2$)
 & $\V_{\vp_2}$ & $n-1$ & $(1+t^5)(1+t^9)\dots (1+t^{4n-3})$ 
\\ 
${\frak so}_{2n+1}$ ($n\ge 2$)
& $\V_{2\vp_1}$ & $n$ & $(1+t^5)(1+t^9)\dots (1+t^{4n+1})$
\\
$\tri$ & $\V_{4\vp}$ & 1 & $1+t^5$ 
\\ 
${\frak so}_{n}$ ($n\ge 5$)
& $\V_{\vp_1}$ & $1$ & $1+t^n$ 
\\
$\tri\oplus\tri$ 
& $\V_{\vp}\otimes\V_{\vp'}$ & $1$ & $1+t^4$ 
\\ 
${\frak f}_{4}$ & $\V_{\vp_1}$ & $2$ & $(1+t^9)(1+t^{17})$
\\ \hline 
\end{tabular}
}  
\end{table}
\noindent
If $\V$ is a symplectic $\g$-module, then $(\wedge^2\V)^\g\ne 0$.
Thus, $(\wedge^\bullet\V)^\g$ cannot be free unless $\dim\V= 2$. 
If $\V$ is neither orthogonal nor symplectic, then all known instances
of free algebras of skew-invariants are those with 
$\dim(\wedge^\bullet\V)^\g=2$.

\sekt{`Spin' of an orthogonal $\Frak g$-module 
and its properties\label{spin1}}%
Let $\V$ be a $\Bbbk$-vector space endowed with a non-degenerate 
quadratic form $Q$. Denote by ${\frak so}(\V)=
{\frak so}_Q(\V)$ the respective
orthogonal Lie algebra and by $\cC_Q(\V)$ the Clifford algebra of $Q$.
Let $\W,\W'$ be maximal $Q$-isotropic subspaces of $\V$ and $\W\cap\W'=0$.
The following relations
are well-known in the theory of Clifford algebras (see e.g. 
\cite[\S\,20.1]{FH}) :

{\sf (i)} $\cC_Q(\V)\simeq \mbox{End\,}(\wedge^{\bullet} \W)$, 
if $\dim\V$ is even,

{\sf (ii)} $\cC_Q(\V)\simeq \mbox{End\,}(\wedge^{\bullet} \W)\oplus
\mbox{End\,}(\wedge^{\bullet} \W')$, 
if $\dim\V$ is odd.
\\[.5ex]
As $\wedge^{\bullet} \W$ (or $\wedge^\bullet\W'$)
is the underlying space of the spin representation of
${\frak so}(\V)$ (in case (i) this representation is the sum of two half-spin
representations), we shall write $Spin(\V)$ in place of
$\wedge^{\bullet} \W$. The above relations are thought of as isomorphisms of 
${\frak so}(\V)$-modules. It is well-known (and easily seen) that $\cC_Q(\V)$
has an ${\frak so}(\V)$-stable filtration such that the associated graded 
algebra is isomorphic to
the exterior algebra of $\V$. Since in both cases $Spin(\V)$ is a self-dual
module, we obtain the following isomorphisms of ${\frak so}(\V)$-modules:
\beq  \label{exterior}
\begin{array}{l}
\wedge^{\bullet} \V\simeq Spin(\V)\otimes Spin(\V), \mbox{ if $\dim\V$ is even}, 
\\
\wedge^{\bullet} \V\simeq 2(Spin(\V)\otimes Spin(\V)), \mbox{ if $\dim\V$ is odd.}
\end{array}
\eeq
Let $\g$ be a reductive Lie algebra and $\pi:\g\rar {\frak so}(\V)$ an
orthogonal representation. Using $\pi$, one may regard $Spin(\V)$ as
$\g$-module. In this way, we obtain a mapping from the set of
orthogonal $\g$-modules to a set of $\g$-modules: $\V\mapsto Spin(\V)$. 
Of course, the $\g$-modules
of the form $Spin(\V)$ must satisfy some constraints; e.g. $\dim Spin(\V)$
is a power of 2. Equations~\re{exterior}, which can be treated as isomorphisms
of $\g$-modules, suggest that `$Spin$' could be used for better
understanding of $\g$-module structure of the exterior algebra of an 
orthogonal module. \\
The point of departure for our considerations is a simple 
formula for the character of
the $\g$-module $Spin(\V)$. Fix some notation, which applies to arbitrary
$\g$-modules (i.e. not necessarily orthogonal ones).
Let ${\cal P}(\V)$ (resp. $\Delta (\V)$) denote the 
set of all (resp. all nonzero) weights of $\V$. 
For instance, $\Delta(\g)=\Delta$.
For $\mu\in {\cal P}(\V)$, $\V^\mu$ is the
corresponding weight space and $m(\mu)=\dim\V^\mu$. If $\V=\V_\lb$ is 
irreducible, 
then the multiplicity is denoted by
$m_\lb(\mu)$. Recall that $\V$ is self-dual if and only if 
$\Delta(\V)=-\Delta(\V)$ and 
$m(\mu)=m(-\mu)$ for all $\mu\in\Delta(\V)$. 
\\
Given an orthogonal $\g$-module $\V$,
let $\Delta(\V)^+$ denote an arbitrary subset such that
$\Delta(\V)=\Delta(\V)^+\sqcup (-\Delta(\V)^+)$.
\begin{s}{Lemma}  \label{char}
$\displaystyle 
\ch Spin(\V)=2^{[m(0)/2]}\prod_{\mu\in\Delta(\V)^+} 
(e^{\mu/2}+e^{-\mu/2})^{m(\mu)}$.
\end{s}\begin{proof}
Using \re{ext-char}, one obtains
$ \displaystyle
\ch(\wedge^\bullet\V)=\prod_{\mu\in{\cal P}(\V)}(1+e^\mu)^{m(\mu)}=$
\[
2^{m(0)}\prod_{\mu\in\Delta(\V)^+}[(1+e^\mu)(1+e^{-\mu})]^{m(\mu)}=
2^{m(0)}
\prod_{\mu\in\Delta(\V)^+}(e^{\mu/2}+e^{-\mu/2})^{2m(\mu)} \ .
\]
Since $\dim\V-m(0)$ is even, comparing with Eq.~\re{exterior} 
completes the proof.
\end{proof}%
Roughly speaking, Eq.~\re{exterior} asserts that a ``square root" of 
$\wedge^\bullet\V$ is again a $\g$-module whenever $\V$ is orthogonal. 
Lemma~\ref{char} gives a precise form for this. 
Notice that the transformation from the proof of Lemma 
can be performed for any self-dual $\g$-module $\V$. But the respective
``square root'' does not yield in general the character of a $\g$-module.
\\[.6ex]
It is convenient to omit the numerical factor in $\ch Spin(\V)$.
The remaining expression is still the character of a $\g$-module. This module
is said to be the {\it reduced Spin} of $\V$  and we write $Spin_0(\V)$ for it:
\beq \label{spin0}
\ch Spin_0(\V)=\prod_{\mu\in\Delta(\V)^+}(e^{\mu/2}+e^{-\mu/2})^{m(\mu)} \ .
\eeq
Several easy properties of $Spin_0$ are summarized below.
\begin{s}{Proposition} \label{easy}
Let $\V=\V^{(1)}$ and $\V^{(2)}$ be orthogonal $\g$-modules. Then \par
(i) $\dim Spin_0(\V)=2^{(\dim\V-m(0))/2}$; \par
(ii) $\wedge^\bullet\V\simeq 2^{m(0)}\cdot Spin_0(\V)^{\otimes 2}$; \par
(iii) $Spin(\V)=Spin_0(\V)$ if and only if $m(0)\le 1$; \par
(iv) $Spin_0(\V^{(1)}\oplus \V^{(2)})\simeq Spin_0(\V^{(1)})\otimes
Spin_0(\V^{(2)})$; \par
(v) $Spin_0(\V)$ is a self-dual $\g$-module.

\end{s}\begin{proof*}
This immediately follows from \re{exterior},\,\re{char}, and \re{spin0}.
\end{proof*}%
\begin{rem}{Examples} \label{ex-adj}
{\sf 1}. Our consideration of $Spin(\V)$ was motivated by the following
observation of Kostant, see \cite[p.~358]{cohom} and \cite{cliff}. \\
Suppose $\V=\g$ and $\pi=\ad$ is the adjoint representation. Then 
$\wedge^\bullet \g\simeq 2^{{\mathrm{\scriptsize rk\,}}\g}(\V_\rho
\otimes \V_\rho)$. 
This means that $Spin(\g)=2^{[{\mathrm{\scriptsize rk\,}}\g/2]}\,\V_\rho$
and $Spin_0(\g)=\V_\rho$. 
\par
{\sf 2}. $\g={\frak sl}_2$. We shall write $R_d$ in place of
$\V_{d\vp_1}$. Recall that $R_2=\g$, $R_d={\cal S}^d R_1$, and $R_d$ is 
orthogonal 
if and only if $d$ is even. Let ${\cal P}(R_1)=\{\esi,-\esi\}$.
Applying Prop.~\ref{easy}, we obtain $Spin_0 (R_{2d})=Spin\, R_{2d}$ and
$\ch Spin\, R_{2d}=\prod_{k=1}^d(e^{k\esi}+e^{-k\esi})$.
It is not hard to compute this character for small values of $d$. Here are
first few formulas: 
$Spin\, R_2{=}R_1$, $Spin\, R_4{=}R_3$, $Spin\, R_6{=}R_6{+}R_0$,
$Spin\, R_8{=}R_{10}{+}R_4$, $Spin\, R_{10}{=}R_{15}{+}R_9{+}R_5$.
It is easily seen that if $Spin\,R_{2d}=\oplus_{i\in I}R_{m_i}$, then
$Spin\,R_{2(d+1)}\supset\oplus_{i\in I}R_{m_i+d+1}$. Therefore the number
of summands is a nondecreasing function of $d$.
\par
{\sf 3}. Let $\W$ be an arbitrary $\g$-module. We may regard $\V:=\W\oplus \W^*$ as 
orthogonal $\g$-module equipped with the quadratic form $Q((w,w^*)):=
\langle w\vert w^*\rangle$, where $(w,w^*)\in \V$ and $\langle\ \vert\ \rangle$
is the canonical pairing of $\W$ and $\W^*$. Assuming for simplicity that
the weights of $\W$ and $\W^*$ are distinct, we see that ${\cal P}(\W)$ can
be taken as $\Delta(\V)^+$. Therefore \par
$\displaystyle
\ch Spin(\V)=\prod_{\mu\in {\cal P}(\W)}(e^{\mu/2}+e^{-\mu/2})^{m(\mu)}=
e^{-\nu}\prod_{\mu\in {\cal P}(\W)}(1+e^{\mu})^{m(\mu)}$, where
$\nu=\frac{1}{2}\sum_{\mu\in {\cal P}(\W)}m(\mu)\mu$. Whence
\[
Spin(\W\oplus \W^*)\simeq \Bbbk_{-\nu}\otimes\wedge^\bullet \W
\simeq  \Bbbk_{\nu}\otimes\wedge^\bullet \W^* \ ,
\]
where $\Bbbk_{-\nu}$ is 1-dimensional $\g$-module with character $-\nu$.
Obviously, $\nu=0$ if and only if $\g\subset {\frak sl}(\W)$, e.g. $\g$
is semisimple. It is not hard to 
verify that the above formula for $Spin(\W\oplus \W^*)$ remains true for all
$\W$.
\end{rem}%
{\bf Definition.}
An orthogonal $\g$-module $\V$ is said to be {\it co-primary}, if
$Spin_0(\V)$ is irreducible. 
\\[1ex]
In this case $Spin(\V)$ is a primary $\g$-module. 
We are going to list all co-primary modules for the semisimple Lie algebras.
At the moment, the following examples of such modules are known:
$\V=\g$, $\g$ simple; $\V=R_4$, $\g={\frak sl}_2$. As a step towards a
classification, we describe another series of co-primary modules. 
\\[.6ex]
Let $\g$ be a simple Lie algebra having two root lengths. We use subscripts
`s` and `l` to mark objects related to short and long roots, respectively.
For instance, $\Delta_s$ is the set of short roots, $\Delta=\Delta_s\sqcup
\Delta_l$, and $\Pi_s=\Pi\cap \Delta_s$. Set $\rho_s=\frac{1}{2}\vert
\Delta^+_s\vert$ and $\rho_l=\frac{1}{2}\vert
\Delta^+_l\vert$. As usual, $s_\ap\in W$ is the reflection corresponding
to $\ap\in\Delta$ and $s_i:=s_{\ap_i}$.
\begin{s}{Lemma}   \label{rhos}
$\rho_s=\sum_{\ap_i\in \Pi_s}\vp_i$.
\end{s}\begin{proof}
It is easily seen that $s_i(\rho_s)=\left\{\begin{array}{cl}
\rho_s , & \mbox{if } \ap_i\in\Pi_l \\
\rho_s-\ap_i , & \mbox{if } \ap_i\in\Pi_s \end{array}\right. $.
\end{proof}%
Let $\theta\in\Delta^+$ be the highest root and $\theta_s$ the short dominant root.
Recall that $\Delta_l=W{\cdot}\theta$, $\Delta_s=W{\cdot}\theta_s$, and 
${|\!|\theta|\!|^2}/{|\!|\theta_s|\!|^2}=2$ or 3. If $\mu\in\Delta$, then
$\mu^\vee:=2\mu/|\!|\mu |\!|^2$.
\begin{s}{Lemma} \label{even}
Suppose ${|\!|\theta|\!|^2}/{|\!|\theta_s|\!|^2}=2$ and $\mu\in \Delta_s$. Then
$(\rho+\rho_s,\mu^\vee)$ is even.
\end{s}\begin{proof}
Let $\mu=\sum_{\ap_i\in\Pi_s}n_i \ap_i+\sum_{\ap_j\in\Pi_l}m_j \ap_j$.
Then $\mu^\vee=
\sum_{\ap_i\in\Pi_s}n_i \ap_i^\vee+2\sum_{\ap_j\in\Pi_l}m_j \ap_j^\vee$.
Therefore $(\rho+\rho_s,\mu^\vee)=(2\rho_s+\rho_l,\mu^\vee)=
(2\rho_s, \sum_{\ap_i\in\Pi_s}n_i \ap_i^\vee)+
(\rho_l, 2\sum_{\ap_j\in\Pi_l}m_j \ap_j^\vee)=
2(\sum_i n_i +\sum_j m_j)$.
\end{proof}%
The following assertion can be proved using classification, but we give
a unified proof.
\begin{s}{Proposition} \label{shortchar}

{\sf (i)} $\dim\V_{\theta_s}=(h+1)m_{\theta_s}(0)$, where $h$ is the 
Coxeter number of $\g$;

{\sf (ii)} $m_{\theta_s}(0)=\# \Pi_s$.
\end{s}\begin{proof} (i)
It is clear that ${\cal P}(\V_{\theta_s})=\{0\}\cup\Delta_s$. Moreover,
$m_{\theta_s}(\ap)=1$ for all $\ap\in\Delta_s$. Applying Freudenthal's 
multiplicity formula
to $m_{\theta_s}(0)$, we obtain  
\[
(\theta_s+2\rho,\theta_s)m_{\theta_s}(0)=2\sum_{\ap\in\Delta^+}\sum_{t\ge 1}
m_{\theta_s}(t\ap)(t\ap,\ap)=2\sum_{\ap\in\Delta_s^+}m_{\theta_s}(\ap)(\ap,\ap)
=2\sum_{\ap\in\Delta_s^+}(\ap,\ap) \ .
\]
Whence 
\[
(1+(\rho,\theta_s^\vee))m_{\theta_s}(0)=\#\Delta_s=\dim\V_{\theta_s}-m_{\theta_s}(0) \ .
\]
As $\theta_s^\vee$ is the highest root in the dual root system
$\Delta^\vee$, we have
$(\rho,\theta_s^\vee)=h-1$.
\par (ii) By part (i), we have $m_{\theta_s}(0)=
\displaystyle \frac{\dim\V_{\theta_s}-m_{\theta_s}(0)}{h}=\frac{\#\Delta_s}{h}$.
Let $c\in W$ be a Coxeter element associated with $\Pi$. It is known
that each orbit of $c$ in $\Delta$ has cardinality $h$ and contains a
unique simple root, see \cite[ch.VI,\,\S\,1,\,Prop.\,33]{4-6}. 
Hence $\#\Delta_s=h(\#\Pi_s)$.
\end{proof}
Some authors call $\V_{\theta_s}$ the {\it little adjoint module\/}. To
a great extent, properties of $\V_{\theta_s}$ are similar with properties
of $\g$.
\begin{s}{Theorem}   \label{short}
Suppose ${|\!|\theta|\!|^2}/{|\!|\theta_s|\!|^2}=2$. Then
$\wedge^\bullet \V_{\theta_s}\simeq 2^{\#\Pi_s}\cdot 
(\V_{\rho_s}\otimes \V_{\rho_s})$.
\end{s}\begin{proof*} 
By Proposition~\ref{shortchar}, we have
$\ch \V_{\theta_s}=\#\Pi_s+\sum_{\ap\in\Delta_s}e^\ap$. Therefore
\[
\ch\wedge^\bullet\V_{\theta_s}=2^{\#\Pi_s}\prod_{\ap\in\Delta_s}(1+e^\ap)=
2^{\#\Pi_s}\prod_{\ap\in\Delta_s^+}(e^{\ap/2}+e^{-\ap/2})^2 \ .
\]
Thus, the statement of theorem is equivalent to that
\beq \label{charros}
\ch\V_{\rho_s}=\prod_{\ap\in\Delta_s^+}(e^{\ap/2}+e^{-\ap/2})=
e^{\rho_s}\prod_{\ap\in\Delta_s^+}(1+e^{-\ap}) .
\eeq
By Weyl's character formula
\beq \label{weyl}
\ch\V_{\rho_s}=\frac{\sum_{w\in W}\esi(w)e^{w(\rho+\rho_s)}
}{\sum_{w\in W}\esi(w)e^{w\rho}}=
\frac{\sum_{w\in W}\esi(w)e^{w(\rho+\rho_s)}}{
\prod_{\ap\in\Delta^+}(e^{\ap/2}-e^{-\ap/2})} \ .
\eeq
Here $\esi(w)=(-1)^{l(w)}$, where $l(w)$ is the length of $w$ with respect
to $\Delta^+$.
Take $\ap\in\Delta_s^+$. We are going to prove that $1+e^{-\ap}$ divides
$\ch\V_{\rho_s}$ in ${\Bbb Z}[{\cal P}]$. Since $\ch\V_{\rho_s}$ is 
$W$-invariant, it is enough to consider the case in which $\ap$ is simple,
i.e., $\ap\in\Pi_s$. Actually, we shall prove
that $1+e^{-\ap}$ divides the numerator in Eq.\,\re{weyl}. For this, we show
how to group together the summands of 
the numerator. 
Let $W^\ap=\{w\in W\mid w^{-1}\ap\in\Delta^+\}$. Then $W$ is the disjoint
union of
pairs $\{s_\ap w,w\}$ ($w\in W^\ap$). Consider the corresponding pairs of 
summands in the numerator of \re{weyl}. Since $\ap\in\Pi_s$, we have 
$\esi(s_\ap w)=-\esi(w)$ and
\[
\esi(w)e^{w(\rho+\rho_s)}+\esi(s_\ap w)e^{s_\ap w(\rho+\rho_s)}=
\esi(w)e^{w(\rho+\rho_s)}(1-e^{-n\ap}) \ ,
\]
where $n=( w(\rho+\rho_s),\ap^\vee)=
( \rho+\rho_s,(w^{-1}\ap)^\vee)$. By the definition of $W^\ap$, 
$n$ is positive. The divisibility will follow
from the fact that $n$ is even. But this is just Lemma~\ref{even}.
\\[.7ex]
Since ${\Bbb Z}[{\cal P}]$ is factorial and the factors $1+e^{-\ap}$ 
($\ap\in\Delta_s^+$) are coprime (see \cite[ch.\,VI,\,\S\,3,\,Lemma~1]{4-6}), 
$e^{\rho_s}\prod_{\ap\in\Delta_s^+}(1+e^{-\ap})$ divides $\ch\V_{\rho_s}$.
The quotient is a $W$-invariant element of ${\Bbb Z}[{\cal P}]$. 
Comparing the maximal terms in both expressions, we see that the quotient
must be equal to 1.
\end{proof*}
\begin{s}{Corollary}
1. $Spin_0(\V_{\theta_s})=\V_{\rho_s}$; \\ \hspace*{3.6cm}
2. $\dim (\wedge^\bullet \V_{\theta_s})^\g=2^{\#\Pi_s}$.
\end{s}%
\vskip-1.2ex
\begin{rem}{Examples} \label{scope} 
To realize the scope of Theorem~\ref{short}, 
we look at all simple Lie algebras with two root lengths. 
\end{rem}%
1. $\g={\frak sp}_{2n}$. Here $\theta=2\vp_1$, $\theta_s=
\vp_2$, and $\rho_s=\vp_1+\ldots+\vp_{n-1}$. Thus
\[
\wedge^\bullet\V_{\vp_2}=2^{n-1}{\cdot}
(\V_{\vp_1+\ldots+\vp_{n-1}})^{\otimes 2}
\quad\mbox{and}\quad
Spin_0(\V_{\vp_2})=\V_{\vp_1+\ldots+\vp_{n-1}} \ .
\]
2. $\g={\frak f}_{4}$. Here $\theta=\vp_4$, $\theta_s=
\vp_1$, and $\rho_s=\vp_1+\vp_{2}$. Thus
\[
\wedge^\bullet\V_{\vp_1}=4(\V_{\vp_1+\vp_2})^{\otimes 2}
\quad\mbox{and}\quad
Spin_0(\V_{\vp_1})=\V_{\vp_1+\vp_2} \ .
\]
3. $\g={\frak so}_{2n+1}$. Here $\theta=\vp_2$, $\theta_s=
\vp_1$, and $\rho_s=\vp_{n}$. In this case
$Spin_0(\V_{\vp_1})=Spin(\V_{\vp_1})=\V_{\vp_{n}}$ and the formula of
Theorem~\ref{short} is nothing but the second equality in Eq.\,\re{exterior}.
Hence the theorem also yields another approach to defining `Spin'
of an orthogonal representation.
\\[.5ex]
4. $\g=\g_2$. Here $|\!|\theta|\!|^2/|\!|\theta_s|\!|^2=3$ and Theorem \ref{short}
does not apply. In this case $\rho_s=\theta_s=\vp_1$ and $\theta=\vp_2$.
An explicit (easy) computation with characters shows that \par
$\wedge^\bullet\V_{\vp_1}=2{\cdot}(\V_{\vp_1}\oplus\odin)^{\otimes 2}$, i.e.,
$Spin(\V_{\vp_1})=\V_{\vp_1}\oplus\odin$. Hence $\V_{\theta_s}$ is not
co-primary. \\
Here $\odin$ stands for the trivial 1-dimensional module.

\begin{subs}{Another proof of Theorem~\ref{short}}\end{subs}
Making use of Weyl's character formula,
we interpret Eq.\,\re{charros} as Weyl's denominator identity for the 
dual root system.
\\
Recall that $\Delta=\Delta_l\sqcup\Delta_s$ and we assume that 
$|\!|\theta|\!|^2/|\!|\theta_s|\!|^2=2$. The dual root system is therefore
isomorphic to $\widetilde\Delta:=\Delta_l\sqcup 2\Delta_s$.
Here $(\widetilde\Delta)_l=2\Delta_s$ and $(\widetilde\Delta)_s=\Delta_l$.
Since $\widetilde W\simeq W$, Weyl's denominator identity for
$\widetilde\Delta$ reads
\[
\sum_{w\in W}\esi(w)e^{w\tilde\rho}=\prod_{\ap\in\tilde\Delta^+}
(e^{\ap/2}-e^{-\ap/2}) \ .
\]
We have $\tilde\rho=\rho+\rho_s$ on the left hand side and
\[
\prod_{\ap\in\Delta_l^+}(e^{\ap/2}-e^{-\ap/2}){\cdot}
\prod_{\mu\in\Delta_s^+}(e^{\mu}-e^{-\mu})=
\prod_{\ap\in\Delta^+}(e^{\ap/2}-e^{-\ap/2})
\prod_{\ap\in\Delta_s^+}(e^{\mu/2}+e^{-\mu/2})\ 
\]
on the right hand side.
Hence dividing Weyl's identity by
$\prod_{\ap\in\Delta^+}(e^{\ap/2}-e^{-\ap/2})$ yields
$\ch\V_{\rho_s}=\prod_{\mu\in\Delta_s^+}(e^{\mu/2}+e^{-\mu/2})$.
\\[1ex]
{\bf Remark.} The previous argument suggests a proper analogue of
\re{charros} for the exceptional Lie algebra $\g_2$. Here 
the dual root system is isomorphic to 
$\tilde\Delta:=\Delta_l\sqcup 3\Delta_s$ and a
similar transformation proves that
$\ch\V_{2\rho_s}=\prod_{\mu\in\Delta_s^+}(e^\mu+1+e^{-\mu})$.

\sekt{Classification of co-primary $\Frak g$-modules\label{spin2}}%
In this section, $\g$ is a semisimple Lie algebra and $\V$ an orthogonal 
$\g$-module. From Eq.\,\re{spin0} it is clear that $\V^\g$ has no affect
on $Spin_0(\V)$. We may therefore assume that $\V^\g=0$.
\begin{s}{Proposition} \label{irrep}
Suppose $\V$ is co-primary. Then there exist decompositions
$\g=\g_1\oplus\ldots\oplus\g_s$, $\V=\V_1\oplus\ldots\oplus\V_s$ such that 

{\sf (i)} Each $\g_i$ is a (semisimple) ideal of $\g$,

{\sf (ii)} $\g_i$ acts trivially on $\V_j$ ($i\ne j$),

{\sf (iii)} $\V_i$ is an irreducible orthogonal co-primary $\g_i$-module.
\end{s}\begin{proof}
Assume that $\V=\V_1\oplus\V_2$, where $\V_1$ and $\V_2$ are orthogonal
$\g$-modules. It follows from the assumptions and Proposition 
\ref{easy}(iv) that the $\g$-module
$Spin_0(\V_1)\otimes Spin_0(\V_2)$ is irreducible. Since both factors
are non-trivial, the only possibility for this is that $\g=\g_1\oplus
\g_2$, where  $\g_i$ acts trivially on $\V_j$ ($i\ne j$) and
$\V_i$ is a co-primary $\g_i$-module ($i=1,2$). Repeating this procedure,
we obtain a decomposition satisfying (i) and (ii), where each $\V_i$
is orthogonal
co-primary and is not a sum of two proper orthogonal $\g_i$-submodules.
Then either $\V_i$ is irreducible or $\V_i=\W_i\oplus\W_i^*$, where $\W_i$ 
is already irreducible. In the second case, we have $Spin(\V_i)\simeq
\wedge^\bullet \W_i$ (see Example \ref{ex-adj}(3)). 
It is easily seen that the $\g_i$-module $\wedge^\bullet \W_i$ is
never primary, i.e., $Spin_0(\V_i)$ can not be irreducible here.
\end{proof}%
Whenever $(\g,\V)$ admits a decomposition satisfying conditions (i) and (ii) of
the Proposition, this will be denoted by\quad
$(\g,\V)=(\g_1,\V_1)\oplus\ldots\oplus (\g_s,\V_s)$. \\
Notice that if each $\V_i$ is irreducible, then all the summands in the
above decomposition are uniquely determined.
\begin{s}{Lemma} \label{m0}
If $\V_\lb$ is an irreducible co-primary $\g$-module, then $m_\lb(0)\ne 0$.
\end{s}\begin{proof}
If $m_\lb(0)=0$, then $\wedge^\bullet\V_\lb\simeq Spin_0(\V_\lb)^{\otimes 2}$,
see \ref{easy}(ii). Since $\dim (\wedge^\bullet\V_\lb)^\g\ge 2$,
the Schur lemma shows that $Spin_0(\V_\lb)$ cannot be irreducible.
\end{proof}%
It follows from the above two assertions that ${\cal P}(\V)$ lies in the root
lattice whenever $\V$ is co-primary. 
\\[.7ex]
Let us present an explicit way for finding some irreducible constituents 
of $Spin_0(\V)$. To write an expression for $\ch Spin_0(\V)$ in \re{spin0},
we exploited
an arbitrary `half' $\Delta(\V)^+$
of $\Delta(\V)$. However a clever choice of  $\Delta(\V)^+$
will provide us with a maximal term in $\ch Spin_0(\V)$ and hence 
with a highest weight. Take $\nu\in {\cal P}_+$ such that $(\nu,\mu)\ne 0$
for all $\mu\in \Delta(\V)$.  Put
$\Delta(\V)_\nu^+=\{\mu\in\Delta(\V)\mid (\mu,\nu)>0\}$. 
A subset of such form is said to be a {\it dominant half\/} of $\Delta(\V)$.
Set $\Lb_\nu:=\frac{1}{2}\sum_\mu m(\mu)\mu$, where $\mu$ ranges over
$\Delta(\V)_\nu^+$.
\begin{s}{Lemma} \label{highest}
$\Lb_\nu$ is a highest weight of $Spin_0(\V)$.
\end{s}\begin{proof}
We show that $\Lb_\nu$ is dominant and it is a maximal element in
${\cal P}(Spin_0(\V))$. Note that the first part is not tautological.
We exploit formula \ref{spin0} with $\Delta(\V)^+_\nu$:
\[\ch Spin_0(\V)=
\prod_{\mu\in\Delta(\V)^+_\nu}(e^{\mu/2}+e^{-\mu/2})^{m(\mu)} \ .
\]
This shows that $e^{\Lb_\nu}$ occurs in $\ch Spin_0(\V)$ with
coefficient 1, $\displaystyle
(\nu,\Lb_\nu)=\max_{\mu\in {\cal P}(Spin_0(\V))}
(\nu,\mu)$, and $\Lb_\nu$ is the unique element of ${\cal P}(Spin_0(\V))$,
where the maximal value is attained. Let $\Lb'_\nu$ be the dominant
representative in $W{\cdot}\Lb_\nu$. Then $\Lb'_\nu\in
{\cal P}(Spin_0(\V))$ and $\Lb'_\nu-\Lb_\nu=\sum_{\ap_i\in\Pi}n_i\ap_i$
with $n_i\ge 0$. Therefore
$(\nu,\Lb'_\nu)\ge(\nu,\Lb_\nu)$ and hence $\Lb'_\nu=\Lb_\nu$.
\end{proof}%
\refstepcounter{equation}%
{\bf (\theequation)} \label{extreme}
The highest weights of $Spin_0(\V)$ of the form $\Lb_\nu$ are said 
to be {\it extreme\/}.
It is easy to describe all dominant halfs of $\Delta(\V)$ and hence
all extreme weights of $Spin_0(\V)$. Consider the
Weyl chamber $C:={\Bbb Q}_+{\cal P}_+\subset {\cal P}_{\Bbb Q}$ 
and its interior $C^o$.
Let $H_\mu$ denote the hyperplane in ${\cal P}_{\Bbb Q}$
orthogonal to $\mu\in \cal P$. Recall that $C^o$ is the connected 
component\footnote{Strictly speaking, use of the term ``connected component''
is correct only for the {\sl real\/} vector space 
${{\cal P}}_{{\bbb R}}$.}
of ${\cal P}_{\Bbb Q}\setminus\cup_{\gamma\in\Delta}H_\gamma$, containing 
dominant weights. Then the hyperplanes $H_\mu$ ($\mu\in\Delta(\V)$) cut
$C$ in smaller chambers. When $\nu$ varies inside of such a `small' chamber
the corresponding extreme weight does not change. We thus obtain a
bijection
\[
\{ \mbox{extreme weights of }Spin(\V)\}\leftrightarrow
\{\mbox{connected components of }C^o\setminus \bigcup_{\mu\in\Delta(\V)}H_\mu\}
\ .\]
In particular, $Spin_0(\V)$ has a unique extreme weight if and only if
  $\Delta(\V)$ has a unique dominant half if and only if
none of the hyperplanes $H_\mu$ cuts $C^o$.
\begin{s}{Lemma}  \label{proport}
Suppose $\Delta(\V)$ lies in the root lattice. Then: \par
\qquad none of the hyperplanes $H_\mu$ ($\mu\in\Delta(\V)$) cuts $C^o$ 
$\Longleftrightarrow$
$\Delta(\V)\subset\cup_{\ap\in\Delta}{\Bbb Z}\ap$. 
\end{s}\begin{proof*}
``$\Leftarrow$'' This is obvious.\\
``$\Rightarrow$'' Assume that $M:=
\Delta(\V)\setminus\cup_{\ap\in\Delta}{\Bbb Z}\ap
\ne\varnothing$. Let $\mu\in M\cap {\cal P}_+$ 
be an element closest to 0. Write $\mu$ as sum of positive roots with
positive integral coefficients
$\mu=\sum_{i=1}^d k_i\gamma_i$ ($\gamma_i\ne\gamma_j$) and so that
$\sum_i k_i$ is minimal over all such presentations. Then $\gamma_i+\gamma_j$
is not a root, i.e., $(\gamma_i,\gamma_j)\ge 0$. Therefore $(\mu,\gamma_1)>0$
and hence $\mu-\gamma_1\in\Delta(\V)$. As $|\!|\mu-\gamma_1|\!|< |\!|\mu|\!|$, we
obtain $\mu-\gamma_1\in\cup_{\ap\in\Delta^+}{\Bbb N}\ap$. Thus, $k_1=1$,
$d=2$ and,
by symmetry, $\mu=\gamma_1+\gamma_2$. Since $(\gamma_1,\gamma_2)\ge 0$ and
$\gamma_1+\gamma_2$ is not a multiple of a root, it is easily seen that
$(\gamma_1,-\gamma_2)$ is a basis of the root system
$\Delta\cap ({\Bbb Q}\gamma_1+{\Bbb Q}\gamma_2)$. Therefore
$(\gamma_1,-\gamma_2)$ is $W$-conjugate to a pair of simple roots
$(\ap_i,\ap_j)$ (see \cite[ch.\,VI,\,\S\,1,Prop.~24]{4-6}).
Thus, $\ap_i-\ap_j\in\Delta(\V)$ and $H_{\ap_i-\ap_j}$ cuts $C^o$.
\end{proof*}%
\begin{s}{Proposition} \label{2.5}
Let $\V$ be a co-primary faithful irreducible $\g$-module. Then 
\par
$\Delta(\V)\subset\cup_{\ap\in\Delta}{\Bbb Z}\ap$ \quad
and \quad $\g$ is simple.
\end{s}\begin{proof}
By Lemma \ref{m0}, $\Delta(\V)$ lies in the root lattice. Therefore the first
claim readily follows from \re{extreme} and Lemma~\ref{proport}.
Assume that $\g=\g_1\oplus\g_2$ is a sum of two ideals. Then
$\V=\V_1\otimes\V_2$, where $\V_i$ is a non-trivial $\g_i$-module. Obviously,
if $\mu_i\in\Delta(\V_i)$ ($i=1,2$), then $\mu_1+\mu_2\in\Delta(\V)$ 
and it is not a multiple of
a root of $\g$.
\end{proof}%
Now, we are ready to state a classification.
\begin{s}{Theorem}  \label{list}
{\sf (i)} Let $\g$ be semisimple and $\V$ a faithful orthogonal
$\g$-module with $\V^\g=0$. Suppose $\V$ is co-primary. Then
\[
(\g,\V)=(\g_1,\V_1)\oplus\ldots\oplus (\g_s,\V_s) \ ,
\]
where each $\g_i$ is simple and $\V_i$ is irreducible and co-primary. 
Each weight of $\V$ is a multiple of a root of $\g$.
\\
{\sf (ii)} If $\g$ is simple and $\V=\V_\lb$ is irreducible and co-primary, then
the pair $(\g,\lb)$ is one of the following:
\par
(a) $\g$ is any and $\lb=\theta$.
\par
(b) $\g\in\{{\frak so}_{2n+1},{\frak sp}_{2n},{\frak f}_{4}\}$ 
and $\lb=\theta_s$.
\par
(c) $\g={\frak so}_{2n+1}$, $\lb=2\theta_s=2\vp_1$ ($n\ge 2$).
\par
(d) $\g={\frak sl}_2$, $\lb=4\vp_1$.
\end{s}\begin{proof*}
(i) By Proposition~\ref{irrep}, such a decomposition with irreducible and
co-primary summands $\V_i$ exists. The other assertions are proved in
Proposition~\ref{2.5}.\\
(ii) By part (i), we have
$\lb\in\{k\theta, k\theta_s \mid k\in {\Bbb N}\}$. \par
{Let $\rk\g=1$}. It follows from Example~\ref{ex-adj}(2) that the only
co-primary $\tri$-modules are $\tri=R_2$ and $R_4$. 
\par
{Let $\rk\g\ge 2$}. Consider the following possibilities. \\
$\bullet$ $\lb=2\theta$. Take $\ap_i\in\Pi$ such that $(\ap_i,\theta)\ne 0$.
Then $2\theta-\ap_i$ is a weight of $\V_{2\theta}$, which is not a multiple of
a root. Thus, $\V_{2\theta}$ is not co-primary.\\
$\bullet$ $\g={\frak sp}_{2n}$ or ${\frak f}_{4}$ and $\lb=2\theta_s$. 
If $\ap_i$ is the unique simple root such that $(\ap_i,\theta_s)\ne 0$, then
$2\theta_s-\ap_i\in\Delta(\V_{2\theta_s})$ is not a multiple of a root.
\\
$\bullet$ $\g={\frak so}_{2n+1}$ and $\lb=3\theta_s=3\vp_1$. Here
$3\vp_1-\ap_1$ is not proportional to a root.
\\
$\bullet$ $\g=\g_2$. We have already shown in Example~\ref{scope}(4)
that $\V_{\theta_s}$ is not co-primary. 
\\[.7ex]
Obviously, if $\V_{k\lb}$ is not co-primary, then the same holds for
any $m\ge k$. Thus, comparing with results of section~\ref{spin1}, we see that the 
only unclear case is {\sl ii(c)}. Our proof that this module is co-primary
is similar to the first proof of Theorem~\ref{short}. 
It will be given in the next
proposition, where we also compute the reduced $Spin$ of $\V_{2\theta_s}$.
\end{proof*}%
\begin{s}{Proposition} \label{B_n} Let $\g={\frak so}_{2n+1}$. Then
\par
1. $Spin_0(\V_{2\vp_1})=\V_{\rho+2\vp_n}$;\par
2. $\wedge^\bullet \V_{2\vp_1}=2^n {\cdot}(\V_{\rho+2\vp_n})^{\otimes 2}$.
\end{s}\begin{proof}
First, we describe the weight structure of the $\g$-module $\V_{2\vp_1}$.
This is easy, since $\V_{2\vp_1}$ is the Cartan (highest) component in
${\cal S}^2\V_{\vp_1}$. Here 
${\cal P}(\V_{2\vp_1})=\{0\}\cup\Delta\cup 2\Delta_s$.
Hence $\Delta(\V)^+=\Delta^+\cup 2\Delta_s^+$. The non-zero weights are of
multiplicity 1, and $m_{2\vp_1}(0)=n$. Therefore, making use of Eq.~\re{spin0},
we obtain
\[
\ch Spin_0(\V_{2\vp_1}){=}\prod_{\ap\in\Delta^+}(e^{\ap/2}+e^{-\ap/2})
\prod_{\ap\in\Delta_s^+}(e^{\ap}+e^{-\ap}){=}e^{\rho+2\vp_n}
\prod_{\ap\in\Delta^+}(1+e^{-\ap})
\prod_{\ap\in\Delta_s^+}(1+e^{-2\ap})=
\]
\[
=e^{\rho+2\vp_n}
\prod_{\ap\in\Delta_l^+}(1+e^{-\ap})
\prod_{\ap\in\Delta_s^+}(1+e^{-\ap}+e^{-2\ap}+e^{-3\ap})\ .
\]
On the other hand,
\beq \label{weyl2}
\ch\V_{\rho+2\vp_n}=
\frac{\sum_{w\in W}\esi(w)e^{w(\rho+\rho+2\vp_n)}}{
\prod_{\ap\in\Delta^+}(e^{\ap/2}-e^{-\ap/2})} \ .
\eeq
Since $\ch Spin_0(\V_{2\vp_1})$ and $\ch\V_{\rho+2\vp_n}$ have the same 
maximal term $e^{\rho+2\vp_n}$, it suffices to prove that each factor in
the last expression for $\ch Spin_0(\V_{2\vp_1})$ divides 
$\ch\V_{\rho+2\vp_n}$, i.e., the numerator in Eq.~\re{weyl2}.

The same procedure, as in the proof of Theorem~\ref{short}, reduces 
the problem to proving that, for any $w\in W^\ap$,
\[  (w(2\rho+2\vp_n),\ap^\vee)\quad\left\{\begin{array}{lcl}
\mbox{is even} & , & \mbox{ if $\ap\in\Pi_l$} \\
\mbox{is divisible by 4} & , & \mbox{ if $\ap\in\Pi_s$}
\end{array}\right. \ .
\]
That is, we need actually to verify that $(w(\rho+\vp_n),\ap^\vee)$
is even whenever $\ap$ is short. As $\vp_n=\rho_s$ for our $\g$,
this is just Lemma~\ref{even}. 

2. This is a formal consequence of part 1, see Proposition~\ref{easy}.
\end{proof}%
Having obtained the list of all irreducible co-primary modules in 
Theorem~\ref{list}(ii),
it is worth looking it through again in order to find out common
features and latent regularities for the representations in question.
\\[.6ex]
First, item {\sl (ii)d\/} in \re{list} can be thought of as starting
point for the series in {\sl (ii)c}. Indeed, $V_{2\vp_1}$ is the
Cartan component in ${\cal S}^2\V_{\vp_1}$ and  $\V_{\vp_1}$ is the
tautological module for ${\frak so}_{2n+1}$ ($n\ge 2$), whereas
$\tri$-module $R_4$ is the Cartan component in ${\cal S}^2R_2$ and
$R_2$ is the tautological module for ${\frak so}_3$.
Thus, the list consists of three groups of representations: 

1. $(\g,\V_\theta=\g)$;

2. $(\g,\V_{\theta_s})$, where $\g$ is of type $\bf B$, $\bf C$, or $\bf F$;

3. $\g={\frak so}(\W)$ and $\V={\cal S}^2_0(\W)$, where $\dim\W=3,5,7,\ldots$
\\[1ex]
The second (more interesting) observation is that, for all items
$(\g,\V)$ in the list, $\g\rar{\frak so}(\V)$ is the
isotropy representation of an {\sl irreducible\/} symmetric space. 
In other words, $\tilde\g:=\g\oplus\V$ has a structure of irreducible
${\Bbb Z}_2$-graded semisimple Lie algebra. 
More precisely, $\tilde\g$ is simple for items 2 and 3,
and $\tilde\g\simeq\g\oplus\g$ for item 1.
Furthermore, it follows from the well-known classification of symmetric
spaces that items 1--3 correspond exactly to the cases, where
$\g$ is {\it non-homologous to zero\/}\footnote{this 
means that the canonical map of homology spaces $H_*(\g)\rar
H_*(\tilde\g)$ is injective.} in $\tilde\g$. The class of homogeneous
spaces $\tilde G/G$ (not necessarily symmetric ones)
such that $\tilde G, G$ are connected and
$\g$ is non-homologous to zero in $\tilde\g$
 has many nice descriptions. We refer the reader to
\cite[\S\,13,\,n.2]{al} for a thorough treatment in the context
of homogeneous spaces of compact Lie groups. In the symmetric case, 
yet another characterization is that this happens
if and only if $\g$
is determined by a {\it diagram involutory automorphism\/} of $\tilde\g$. 
An explicit description of the diagram automorphisms of simple Lie
algebras is found in \cite[\S\,7.9,\,7.10]{kac}. 
The third observation is that any irreducible co-primary module occurs in
Table~\ref{final} in section~\ref{grassmann}, i.e., it has a free
algebra of skew-invariants. 

These observations give us some hope that the reduced $Spin$ of the isotropy
representation of an {\sl arbitrary\/} symmetric spaces might have some
interesting properties. This is really the case and we turn to such 
considerations in the following sections.

\sekt{Some auxiliary results\label{spin3}}
In this section, we prove an auxiliary result on Weyl groups and
recall some standard facts on involutions of simple Lie algebras.
\\[1ex]
Let $W$ be the Weyl corresponding to a reduced root system $\Delta$ with a set
of positive roots  $\Delta^+$.
Let $w\mapsto l(w)$ be the length function on $W$ determined by $\Delta^+$.
Recall that $l$ can be defined as 
$l(w)=\#\{\ap\in\Delta^+\mid w(\ap)\in\Delta^-\}$.
Consider an arbitrary  subset $\Delta_0\subset\Delta$ which is a root 
system in its 
own right, but is not necessarily closed in $\Delta$. That is, it is allowed
that $\ap+\beta\in \Delta\setminus\Delta_0$ for some $\ap,\beta\in
\Delta_0$. It is easily seen that such a phenomenon can only occur if
$\Delta$ has roots of different length.
As a sample of such non-closed subset, we mention $\Delta_0=\Delta_s$.
Nevertheless, 
$W_0$, the Weyl group of $\Delta_0$, is always
identified with a subgroup of $W$.
Clearly, $\Delta_0^+:=\Delta_0\cap\Delta^+$ can be taken as  set of 
positive roots for $\Delta_0$. 
\begin{s}{Proposition}  \label{coset}
1. Any coset $wW_0\subset W$ contains a unique representative of minimal 
length.
Denoting by $W^0$ the set of minimal length representatives,  we have
$W^0=\{w\in W\mid w(\Delta^+_0)\subset \Delta^+\}$. \\
2. The mapping $W^0\times W_0\rar W$ ($(w^o,w_o)\mapsto w_o(w^o)^{-1}$) is a
bijection.
\end{s}\begin{proof}
1. Set $W'=\{w\in W\mid w(\Delta^+_0)\subset \Delta^+\}$. Then $W^0\subset
W'$. Indeed, assume that $w\in W^0$ and $w(\beta)\in\Delta^-$ for some
$\beta\in\Delta^+_0$. Then $ws_\beta(\beta)\in\Delta^+$ and it follows from
\cite[2.3]{bgg} that $l(ws_\beta)< l(w)$. But this contradicts the fact that
$w\in W_0$. 
Obviously, each coset contains elements of minimal length and hence 
elements from $W'$.
Assume that $u,v
\in W'\cap vW_0$. Then $u=vw$ for some $w\in W_0$. If $w\ne e$, then
$w(\beta)\in\Delta^-_0$ for some $\beta\in\Delta^+_0$. 
Whence $u(\beta)=v(w(\beta))
\in\Delta^-$, which contradicts the
assumption. Thus, each coset contains a unique element of $W'$, $W^0=W'$,
and we are done. \par
2. Obvious.
\end{proof}%
{\bf Remark.} If $\Delta_0$ is generated by a part of the basis $\Pi\subset
\Delta^+$ 
(i.e., $\Delta_0\cap\Pi$ is a basis of $\Delta_0$), then $W_0$ is a parabolic
subgroup of $W$. In this case the Proposition is well known and, moreover, 
the relation $l(w^ow_0)=l(w^o)+l(w_o)$ holds, see e.g. \cite[1.10]{hump}.
However this relation does not hold in general. 
\\
For $w_o\in W_0$, let $l_0(w_o)$ denote the length of $w_o$ in $W_0$.
That is, $l_0(w_o)=\#\{\mu\in\Delta_0^+\mid w_o(\mu)\in\Delta_0^-\}$. 
If $W_0$ is a parabolic subgroup, then
$l_0(w_o)=l(w_o)$, but in general we have only ``$\le$''. The usual 
{\it determinant\/} or {\it parity\/} for the elements of $W$ is defined by
$\esi(w)=(-1)^{l(w)}$. Making use of the above bijection, one may introduce
a parity depending on $\Delta_0$.
By Prop.~\ref{coset}(2), each element $w\in W$ has a unique presentation
$w=w_o(w^o)^{-1}$, where $w_o\in W_0$ and $w^o\in W^0$. 
Set $l_0(w):={l_0(w_o)}$ and $\tau(w):=(-1)^{l_0(w)}$. 
So, if $w=w_o$, then $\tau(w_o)$ is nothing but the usual parity on $W_0$,
which will be denoted by $\esi_0(w_o)$. Therefore one may say that $\tau$ 
is the extension of the parity $\esi_0$ to $W$ determined by the `section'
$W^0$. 
The function $w\in W\mapsto \tau(w)\in\{1,-1\}$ is said to be the 
{\it cunning parity\/} on $W$, determined by $\Delta_0^+$ (or $W_0$). 
It is convenient
to give an expression for $l_0(w)$, and hence for $\tau(w)$, where
$w_o$ is not explicitly mentioned.
\begin{s}{Lemma}  \label{tau}
$l_0(w)=\#\{\ap\in\Delta^-\mid w(\ap)\in\Delta_0^+\}$.
\end{s}\begin{proof*}
Let $w=w_o(w^o)^{-1}$, as above. Consider the subsets \par 
$M_1=\{\ap\in\Delta^-\mid w(\ap)\in\Delta_0^+\}$ \quad and \quad
$M_2=\{\mu\in\Delta^-\mid w_o(\mu)\in\Delta_0^+\}$. \\
Since $\Delta_0$ is $W_0$-stable, $M_2\subset\Delta_0^-$ and therefore
$l_0(w_o)=\# M_2$. By Prop.~\ref{coset}(1), we have
$(w^o)^{-1}M_1\cap\Delta_0^+=\varnothing$. Since $w_o((w^o)^{-1}M_1)
\subset \Delta^+_0$, we see that $(w_o)^{-1}M_1\subset\Delta_0^-$. Thus
$(w^o)^{-1}M_1\subset M_2$. Similarly, one proves the opposite containment.
Thus, $l_0(w)=\# M_2=\# M_1$, and we are done.
\end{proof*}%

\subs{Classes of involutory automorphisms} \label{classes}
Here $\g$ is a simple Lie algebra.
\\[1ex]
Given an involutory automorphism $\Theta$ of $\g$, consider the 
${\Bbb Z}_2$-grading
$\g=\g_0\oplus\g_1$, where $\g_i=\{x\in\g\mid \Theta(x)=(-1)^ix\}$. 
The reductive subalgebra $\g_0$ is called {\it symmetric\/}. 
The involutory automorphisms
fall into three classes:

a) $\rk\g=\rk\g_0$ and $\g_0$ is semisimple;

b) $\rk\g=\rk\g_0$ and $\g_0$ has 1-dimensional centre;

c) $\rk\g >\rk\g_0$. \\
In cases a) and b), $\Theta$ is inner and, accordingly, both $\g_0$ 
and the ${\Bbb Z}_2$-grading are said to be
{\it of inner type\/}. It is well known that the $\g_0$-module $\g_1$
is irreducible in cases a) and c), and is the sum of two dual 
submodules in case
b). However, $\g_1$ is orthogonal in all three cases and one may
consider the $\g_0$-module $Spin(\g_1)$. An important feature of this
situation is that all nonzero weights of $\g_1$ are of multiplicity 1. 
This is clear in the equal rank cases, and can
also be proved for c). An invariant theoretic proof of this uses 
Lemma~3.4 in \cite{k80} and the fact that
the linear group $G_0\rar GL(\g_1)$ is {\it visible\/}.

\sekt{ ${\bf Spin}({\Frak g}_1)$ for the inner involutory 
automorphisms\label{spin4}}%
In this section, $\g$ is simple and $\g_0$ is a symmetric subalgebra of 
inner type. Retain for $\g$ the previous notation such as $\te$, $\Delta$,
$\Delta^+$, $\rho$, $C$, etc. Since $\rk\g=\rk\g_0$, we may assume 
that $\te$ is a Cartan
subalgebra in both $\g$ and $\g_0$. Let $\Delta_0$ be the root system of
$(\g_0,\te)$ and $\Delta_1$ the set of weights of the $\g_0$-module
$\g_1$. Then $\Delta=\Delta_0\sqcup\Delta_1$ and
$\Delta_0$ is a closed subset of $\Delta$.
We regard $\Delta_0^+:=\Delta^+\cap\Delta_0$ as set of positive roots
for $\g_0$. Note also that $\Delta_1$ contains a distinguished `half'
$\Delta_1^+=\Delta^+\cap\Delta_1$. Then Prop.~\ref{coset} applies to the Weyl
groups $W_0\subset W$ and one obtains the ``minimal length'' 
subset $W^0\subset W$.
\\[.8ex]
Our aim is to describe the $\g_0$-module $Spin_0(\g_1)$. As $\g_1$ has
no zero weight, we have $Spin(\g_1)=Spin_0(\g_1)$.
As a first step,  we find all extreme weights of $Spin(\g_1)$.
Recall from \re{highest},~\re{extreme} that each dominant half of
$\Delta_1$ determines an extreme weight for $Spin(\g_1)$. According
to that discussion, one has to take the dominant chamber $C_0$ for $\g_0$
and cut it up by the hyperplanes orthogonal to the roots of $\Delta_1$.
Clearly, each small chamber is isomorphic to $C$. Since there are $\# W$
chambers for $\g$ and $\# W_0$ chambers for $\g_0$, we obtain the partition of
$C_0$ in $\# (W/W_0)$ small chambers. Then any weight inside of a
small chamber determines a dominant half of $\Delta_1$ and an extreme weight.
In the following proposition we give a formula for these extreme weights.
\begin{s}{Proposition} \label{inn-extr} \\
1. The set of hyperplanes $H_\mu$ ($\mu\in\Delta_1$) cuts $C_0$ in 
$\# (W/W_0)$ small chambers;  \\
2. The collection of ($\g_0$-dominant) weights $w^{-1}\rho$ ($w\in W^0$)
contains representatives of all small chambers in $C_0$. \\
3. The extreme weight of $Spin(\g_1)$ corresponding to $w^{-1}\rho$ is
$\lb_w:=w^{-1}\rho-\rho_0$.
\end{s}\begin{proof}
1. This is proved in the previous paragraph. \\
2 \& 3. If $\ap\in\Delta_0^+$ and $w\in W^0$, then $w\ap\in\Delta^+$,
see Prop.~\ref{coset}(1). Therefore $w^{-1}\rho$ is $\g_0$-dominant.
Since the number of these weights is $\# (W/W_0)$, as required, it suffices
to verify that the corresponding dominant halfs are different.
\\[.7ex]
By definition, the dominant half of $\Delta_1$ associated with $w^{-1}\rho$
is 
\[
(\Delta_1)_{w}^+=\{
\mu\in\Delta_1\mid (w^{-1}\rho,\mu)>0\}=\{
\mu\in\Delta_1\mid w\mu\in\Delta^+ \} \ .
\]
Because all weight multiplicities in $\g_1$ are equal
to 1, the corresponding extreme weight is
$\lb_w:=\frac{1}{2}|(\Delta_1)_{w^{-1}\rho}^+|$. Set
$M_w=\{\mu\in\Delta_1^+\mid w\mu\in\Delta^+\}$ and 
$\ov{M}_w=\Delta_1^+\setminus M_w$. Then $\Delta^+=\Delta_0^+\sqcup
M_w\sqcup \ov{M}_w$ and 
\[
\rho=\rho_0+ \textstyle
\frac{1}{2}|M_w|+\frac{1}{2}|\ov{M}_w| \ .
\]
Since $w\in W^0$, we obtain
\[
w^{-1}\rho=\rho_0+ \textstyle
\frac{1}{2}|M_w|-\frac{1}{2}|\ov{M}_w| \ .
\]
Note also that $|(\Delta_1)_{w}^+|=|M_w|-|\ov{M}_w|$. Whence
$\lb_{w}=w^{-1}\rho-\rho_0$. Thus, we have obtained the required
number of different extreme weights.
\end{proof}%
In the remainder of the section, notation $\V_\lb$ refers to a $\g_0$-module.

\begin{s}{Theorem}  \label{main-inner}
Let $\g=\g_0\oplus\g_1$ be a ${\Bbb Z}_2$-grading of inner type. Then \\[1.8ex]
\centerline{ $\displaystyle
Spin_0(\g_1)=Spin(\g_1)=\bigoplus_{w\in W^0}\V_{\lb_w}$ \ .} \vskip-.6ex 
\end{s}\begin{proof*}
It follows from the preceding exposition that
\[
\bigoplus_{w\in W^0}\V_{\lb_w}\subset Spin_0(\g_1)=Spin(\g_1)\ .
\]
Since $Spin(\g_1)$ is self-dual,  $\dim (Spin(\g_1)^{\otimes 2})^{\g_0}$
is greater than or equal to the number of irreducible summands of $Spin(\g_1)$.
Therefore the desired equality is equivalent to that 
$\dim (Spin(\g_1)^{\otimes 2})^{\g_0}=\# W^0$. Recall the main property
of `Spin' in our situation:
\[
\wedge^\bullet\g_1\simeq Spin(\g_1)^{\otimes 2} \ .
\]
Hence the question about $\g_0$-invariants is being translated in the setting
of exterior algebras. Assuming that $\Bbbk={\Bbb C}$, we can
exploit de Rham cohomology with complex coefficients.
It is well known that $(\wedge^\bullet\g_1)^{\g_0}$
is isomorphic to $H^*(G/G_0)$, the cohomology ring of the symmetric 
space $G/G_0$ \cite[\S 9\, n.11]{al}, and that
$\dim H^*(G/G_0)=\# (W/W_0)$ \cite[\S 13\, n.3]{al}. 
This completes the proof.
\end{proof*}%
\begin{rem}{Example} \label{hermitian}
Let $\Theta$ be a `Hermitian' involutory 
automorphism, i.e.,
$\g_0$ has a 1-dimensional centre and $\g_1\simeq \W\oplus\W^*$, where
$\W$ is a faithful irreducible $\g_0$-module. This is just case
\ref{classes}(b). Then $\Delta(\W)={\cal P}(\W)=\Delta_1^+$ and, according
to Example \ref{ex-adj}(3), $Spin(\g_1)\simeq \Bbbk_{\rho_1}\otimes
\wedge^\bullet \W^*$, where $\rho_1=\frac{1}{2}|\Delta_1^+|$. It then
follows from Theorem~\ref{main-inner} that
\[
\wedge^\bullet \W^*=\Bbbk_{-\rho_1}\otimes Spin(\g_1)=\bigoplus_{w\in W^0}
\V_{w\rho-\rho} \ .
\]
Or, equivalently, $\{\rho-w\rho\mid w\in W^0\}$ is the set of all highest
weights for the $\g_0$-module $\wedge^\bullet \W$. 
This result was obtained by Kostant (see \cite[8.2]{cohom}) as application 
of his results on the cohomology of the nilpotent radical of a parabolic
subalgebra of $\g$. In this situation, $\W$ is the Abelian nilpotent radical
of the parabolic subalgebra $\g_0\oplus \W$.
So, the concept of `Spin' and Theorem~\ref{main-inner} yield another
generalization of Kostant's result.
\end{rem}%
Purists may condemn the above proof of Theorem~\ref{main-inner}, 
since it invokes
the cohomology theory of compact Lie groups over $\Bbb C$. Fortunately,
there exists also a rather simple and purely algebraic proof. We shall
show that the equality in \ref{main-inner} is equivalent to an identity
in ${\Bbb Z}[{\cal P}]$, which is a variation of the Weyl denominator
formula. Recall from section~\ref{spin3} the cunning parity $\tau(w)$ for $w\in W$,
determined by $W_0$.
\begin{s}{Theorem}  \label{identity}
Let $\Delta=\Delta_0\sqcup \Delta_1$ be the partition corresponding to 
a ${\Bbb Z}_2$-grading of $\g$ of inner type.  Then \\[1.5ex]
\centerline{
$ \displaystyle \sum_{w\in W}\tau(w)e^{w\rho}=\prod_{\ap\in \Delta_0^+}
(e^{\ap/2}-e^{-\ap/2})\prod_{\mu\in \Delta_1^+}
(e^{\mu/2}+e^{-\mu/2})$\ .
} \vskip-.6ex
\end{s}\begin{proof*}
The fact that $\Delta_0$ and $\Delta_1$ originate from an inner involutory
automorphism can alternatively be stated as follows:  
\\[.5ex]
\hspace*{5pt} ${\bf (\ast)}$ \qquad
\parbox{350pt}{%
if $\ap\in\Delta_i$, $\beta\in\Delta_j$, and $\ap+\beta\in\Delta$, then
$\ap+\beta\in\Delta_{i+j}$, }
\\[1ex]
where, of course, $i,j\in {\Bbb Z}/2{\Bbb Z}$.
Let ${\cal Q}\subset {\cal P}$ be the root lattice. For ${\Bbb Z}[{\cal Q}]$,
with basis $e^\ap$ ($\ap\in {\cal Q}$), one has a version of
Weyl's denominator identity:
\[
\sum_{w\in W}\esi(w)e^{w\rho-\rho}=\prod_{\ap\in \Delta^+}
(1-e^{-\ap}) \ .
\]
Consider the second copy of ${\Bbb Z}[{\cal Q}]$,
with basis $q^\ap$ ($\ap\in {\cal Q}$), and the equality in
${\Bbb Z}[{\cal Q}]\otimes {\Bbb Z}[{\cal Q}]$:
\beq \label{double}
\sum_{w\in W}\esi(w)q^{w\rho-\rho}e^{w\rho-\rho}=\prod_{\ap\in \Delta^+}
(1-q^{-\ap}e^{-\ap}) \ .
\eeq
Take the specialization of this identity such that
$q^\ap \rar\left\{ \begin{array}{rc}1, & \ap\in\Delta_0 \\
-1, & \ap\in \Delta_1 \end{array}\right.$. It has to be verified that
one obtains a well-defined homomorphism
$({\cal Q}, {+})\rar \{1,-1\}\simeq {\Bbb Z}/2{\Bbb Z}$. In other words, if
$\nu=\sum_{i\in I}\mu_i$ is a sum of roots then the number of summands lying
in $\Delta_1$ should have the same parity for all such presentations.
Indeed, assume that $\sum_{i\in I}\mu_i=\sum_{j\in J}\beta_j$. We argue
by induction on $\# I{+}\# J$. Since 
$(\sum_{i\in I}\mu_i,\sum_{j\in J}\beta_j)>0$, there exist $i_0,j_0$
such that $(\mu_{i_0},\beta_{j_0})>0$. Hence $\mu_{i_0}-\beta_{j_0}$
is a root and
$\sum_{i\in I\setminus \{i_0\}}\mu_i+(\mu_{i_0}-\beta_{j_0})=
\sum_{j\in J\setminus \{j_0\}}\beta_j$. We conclude by applying the
inductive hypothesis to this equality and using ${\bf (\ast)}$.
\\
Thus, the specialization is well-defined and we obtain 
$\displaystyle \prod_{\ap\in \Delta_0^+}
(1-e^{-\ap})\prod_{\mu\in \Delta_1^+}
(1+e^{-\mu})$ at the right hand side of Eq.~\re{double}. It is easily
seen that $w\rho-\rho=-|\Delta(w)|$, where $\Delta(w)=\{
\ap\in\Delta^+\mid w^{-1}\ap\in\Delta^-\}=\Delta^+\cap w(\Delta^-)$.
Therefore $q^{w\rho-\rho}$ specializes to $(-1)^n$, where 
$n=\#(\Delta_1^+\cap w\Delta^-)$.  Recall that $\esi(w)=(-1)^{l(w)}$
and $l(w)=\#(\Delta^+\cap w\Delta^-)$. Thus the resulting sign 
on the left hand side is
$(-1)^{\#(\Delta_0\cap w\Delta^-)}$, which is just $\tau(w)$ by 
Lemma~\ref{tau}. This completes the proof of the theorem.
\end{proof*}%
\begin{subs}{Another proof of theorem~\ref{main-inner}} \label{another-inn}
\end{subs}
By Weyl's character formula for $\g_0$-modules and Prop.~\ref{inn-extr}(3),
\[
\ch \V_{\lb_w}=\frac{\sum_{\tilde w\in W_0}\esi_0(\tilde w)
e^{\tilde w(\rho_0+\lb_w)}}{
\prod_{\ap\in\Delta_0^+}(e^{\ap/2}-e^{-\ap/2}) }=
\frac{\sum_{\tilde w\in W_0}\esi_0(\tilde w)e^{\tilde w{w^{-1}}\rho}}{
\prod_{\ap\in\Delta_0^+}(e^{\ap/2}-e^{-\ap/2}) } \ .
\]
Hence
\[
\ch\Bigl(\bigoplus_{w\in W^0}\V_{\lb_w}\Bigr)=
\frac{\sum_{w\in W^0}\sum_{\tilde w\in W_0}\esi_0(\tilde w)
e^{\tilde w{w^{-1}}\rho}}{
\prod_{\ap\in\Delta_0^+}(e^{\ap/2}-e^{-\ap/2}) } \ .
\]
By the very definition of $\tau(w)$ (see section~\ref{spin3}) and
Prop.~\ref{coset}(2), it follows that the numerator
is equal to $\sum_{w\in W}\tau(w)e^{w\rho}$. Whence, by Theorem~\ref{identity},
\[
\ch\Bigl(\bigoplus_{w\in W^0}\V_{\lb_w}\Bigr)=
\prod_{\mu\in\Delta_1^+}(e^{\mu/2}+e^{-\mu/2})=\ch Spin (\g_1) \ .
\]
\refstepcounter{equation}\label{ex-inner}%
{\bf \theequation\quad Examples.}\quad
{\sf 1}. $\g={\frak so}_{2n+1}$, $\g_0={\frak so}_{2n}$.
Here $\g_1\simeq\V_{\vp_1}$ is the tautological ${\frak so}_{2n}$-module
and $\# (W/W_0)=2$. Let $\{\esi_1,\ldots,\esi_n\}$ be the standard basis
of $\te^*$ so that $\Delta=\{\pm\esi_i\pm\esi_j,\, \pm\esi_i \mid
1\le i,j\le n, i\ne j\}$. Here $\Delta_0^+=
\{\esi_i\pm\esi_j\, (i < j), \esi_i\}$ and $\Delta_0=\Delta_l$.
Then $W^0=\{id, w_n\}$, where $w_n(\esi_i)=\esi_i$ ($i\le n-1$)
and $w_n(\esi_n)=-\esi_n$. Since $\Delta_1=\Delta(\V_{\vp_1})=
\{\pm\esi_1,\ldots,\pm\esi_n\}$, the corresponding dominant halfs are
$(\Delta_1)^+_{id}=\{\esi_1,\ldots,\esi_n\}$  and
$(\Delta_1)^+_{w_n}=\{\esi_1,\ldots,\esi_{n-1},-\esi_n\}$, and the
corresponding extreme weights are $\vp_n$ and $\vp_{n-1}$.
Thus, $Spin(\V_{\vp_1})=\V_{\vp_{n-1}}\oplus\V_{\vp_n}$ and
$\wedge^\bullet\V_{\vp_1}=(\V_{\vp_{n-1}}\oplus\V_{\vp_n})^{\otimes 2}$.
Notice that the last equality is nothing but the first equality in 
Eq.~\re{exterior}.
\\[.6ex]
{\sf 2}. $\g={\frak f}_4$, $\g_0={\frak so}_9$. 
Here $\g_1\simeq\V_{\vp_4}$ and $\# (W/W_0)=3$. 
In the standard notation for ${\frak f}_4$, we have
$\Delta^+=\{\esi_i\pm\esi_j \, (i {<} j),\, \esi_i,\, \frac{1}{2}(\esi_1\pm
\esi_2\pm\esi_3\pm\esi_4)\}$. Then $\Delta_0^+=\Delta_l^+
\sqcup\{\esi_1,\esi_2,\esi_3,\esi_4\}$. An explicit computation shows that
$W^0=\{id,w',w''\}$, where\\[.6ex]
\centerline{
$w' : \left\{
\begin{array}{c}
\esi_1\mapsto \frac{1}{2}(\esi_1+\esi_2+\esi_3+\esi_4) \\
\esi_2\mapsto \frac{1}{2}(\esi_1+\esi_2-\esi_3-\esi_4) \\
\esi_3\mapsto \frac{1}{2}(\esi_1-\esi_2+\esi_3-\esi_4) \\
\esi_4\mapsto \frac{1}{2}(\esi_1-\esi_2-\esi_3+\esi_4) 
\end{array}
\right.$ and 
$w'' : \left\{
\begin{array}{c}
\esi_1\mapsto \frac{1}{2}(\esi_1+\esi_2+\esi_3-\esi_4) \\
\esi_2\mapsto \frac{1}{2}(\esi_1+\esi_2-\esi_3+\esi_4) \\
\esi_3\mapsto \frac{1}{2}(\esi_1-\esi_2+\esi_3+\esi_4) \\
\esi_4\mapsto \frac{1}{2}(\esi_1-\esi_2-\esi_3-\esi_4) 
\end{array}
\right.$.}
\\[.6ex] (One may notice that any $w\in W^0$ must preserve $(\Delta_0)_l^+=
\Delta_l^+$, $\Delta_l$ being the root system of type $\GR{D}{4}$.
Hence $w$ takes $\esi_2-\esi_3$ to itself and 
permutes somehow $\esi_1-\esi_2,\, \esi_3-\esi_4$, and $\esi_3+\esi_4$.)
Whence \par
$\lb_{id}=\rho-\rho_0=2\esi_1$,

$\lb_{w'}=(w')^{-1}\rho-\rho_0=\esi_1+\esi_2+\esi_3$,

$\lb_{w''}=(w'')^{-1}\rho-\rho_0=(3\esi_1+\esi_2+\esi_3+\esi_4)/2$. \\
Thus \\[.6ex]
\centerline{ $Spin(\V_{\vp_4})=\V_{2\vp_1}\oplus\V_{\vp_3}\oplus
\V_{\vp_1+\vp_4}$.}
\vskip2ex

\sekt{ ${\bf Spin}({\Frak g}_1)$ for the outer 
involutory automorphisms\label{spin5}}%
In this section $\g$ is simple and $\Theta$ is outer.
Since $\rk\g_0<\rk\g$, there is no clear relation between roots and Weyl groups
of the two algebras, and the approach of section~\ref{spin4} seems to fail completely.
Yet, it appears to be possible to describe $Spin(\g_1)$ in a similar fashion,
but with some complications. Another price is that we have to exploit
case-by-case arguments several times.

\begin{subs}{Associated diagram involutory automorphism of $\g$}
\end{subs}
By a result of Steinberg, $\Theta$ keeps stable a Borel subalgebra and
a Cartan subalgebra in it. Therefore we may (and shall) assume that
$\Theta\te=\te$ and $\Theta\ut^+=\ut^+$. Then $\Theta$ also preserves
$\Delta^+$ and $\Pi$, as subsets of $\te^*$. In particular, $\Theta$
induces an involution of the Dynkin diagram. Associated with this
involution,
one has the specific involutory automorphism of $\g$, which is called
the {\it diagram involutory automorphism\/} and denoted by $\OT$. 
Roughly speaking, $\OT$ performs the same involution on $\Pi$, as $\Theta$, 
and transforms `well' the Chevalley generators of $\g$\footnote{Explicit 
formulas for the diagram automorphisms of all simple Lie algebras are 
written in \cite[\S\,7.9,7.10]{kac}.}.
We are going to compare properties of the ${\Bbb Z}_2$-gradings 
\[
\g=\g_0\oplus\g_1  \quad\mbox{and}\quad \g=\g_{\ov{0}}\oplus\g_{\ov{1}}
\]
arising from $\Theta$ and $\OT$.
By construction, $\Theta\vert_\te=\OT\vert_\te$.
Therefore $\Theta$ and $\OT$ act identically on $\Delta^+$ and
$\te_0:=\te^\Theta$ is a Cartan subalgebra for both $\g_0$ and $\g_{\ov{0}}$.
Let us organize notation for roots and weights of the symmetric
subalgebras in question:

$\bullet$ \ $\Delta_0$ (resp. $\Delta_{\ov{0}}$) is the root system
of $\g_0$ (resp. $\g_{\ov{0}}$) relative to $\te_0$;

$\bullet$ \ $\Delta_1$ (resp. $\Delta_{\ov{1}}$) is the set of non-zero
weights of the $\g_0$-module $\g_1$
(resp. $\g_{\ov{0}}$-module $\g_{\ov{1}}$) relative to $\te_0$.
\\
Since all these sets are defined with respect to
a common Cartan subalgebra, 
$\Delta_0\cup\Delta_1=\Delta_{\ov{0}}\cup\Delta_{\ov{1}}$ and, more
precisely, the totality of weights occurring in $\{\Delta_0,\Delta_1\}$
is the same as in $\{\Delta_{\ov{0}},\Delta_{\ov{1}}\}$.
Because $\te_0$ contains regular elements of $\g$ 
(see e.g. \cite[8.1(b)]{kac}),
none of the roots of $\g$ vanishes on $\te_0$. Therefore the above
totality of weights consists of all restricted roots.
Moreover, since the non-zero weights of $\te_0$ in $\g_1$ (or $\g_{\ov{1}}$)
are of multiplicity 1, 
\[
\#\Delta=\# \Delta_0+\#\Delta_1=\#\Delta_{\ov{0}}+\#\Delta_{\ov{1}} \ .
\]
{\bf Warning.}
Unlike section~\ref{spin4}, elements of $\Delta_0$ and $\Delta_1$ have not
much in common with roots of $\g$. 
Actually, we do not need $\Delta$ in this section.
\\[.6ex]
The next assertion follows from the classification. 
\\[1ex]
{\bf Fact.} 
The fixed point subalgebra of an outer involutory automorphism always
has roots of different
length, with $|\!|long|\!|^2/|\!|short|\!|^2=2$.
\\[1ex]
This applies to both $\g_0$ and $\g_{\ov{0}}$ and,
as in section~\ref{spin1}, we use the subscripts `s' and `l' to denote the objects
related to short and long roots in $\Delta_0$ and $\Delta_{\ov{0}}$.
A close look to the classification list 
reveals important features of this situation. 
\\[1ex]
\refstepcounter{equation}
{\bf (\theequation)}\ \label{look} Suppose $\Theta\ne\OT$. Then
$\left\{\begin{array}{l}
\mbox{
 $\g_{\ov{0}}$ is simple and $\g_{\ov{1}}$ is the
little adjoint module for $\g_{\ov{0}}$;}\\
\mbox{
 $\Delta_{\ov{1}}=(\Delta_{\ov{0}})_s$;}\\
\mbox{
 $\Delta_{0}\subset\Delta_{\ov{0}}$ and
$(\Delta_{0})_s=(\Delta_{\ov{0}})_s$. }
\end{array}\right.$ \\[.5ex]
In fact, there are 7 series of outer involutory automorphisms of simple
Lie algebras. They form three pairs $(\Theta,\OT)$ and one ``isolated''
diagram involutory automorphism, where \re{look} is not satisfied.
The relevant data for all these series are presented in Table \ref{thety}.

\begin{table}[htbp]
\caption{The outer involutory automorphisms} \label{thety}
\vskip.8ex\centerline
{
\begin{tabular}{c|cc|cc|c}
& \multicolumn{2}{c|}{$\Theta$} & \multicolumn{2}{c|}{$\OT$} & \\ \hline 
$\g$ & $\g_0$ & $\g_1$ & $\g_{\ov{0}}$ & $\g_{\ov{1}}$ 
& $\# W_{\ov{0}}/W_0$ \\ \hline
\hline 
${\frak sl}_{2n}$ & ${\frak so}_{2n}$ & 
$\V_{2\vp_1}$ & ${\frak sp}_{2n}$ & $\V_{\vp_2}$ & 2 \\ 
${\frak s0}_{2n+2m+2}$ & ${\frak so}_{2n+1}\oplus {\frak so}_{2m+1}$ & 
$\V_{\vp_1}\otimes\V'_{\vp_1}$ & ${\frak so}_{2n+2m+1}$ & $\V_{\vp_1}$ & 
${n+m \choose m}$  \\ 
${\frak e}_{6}$ & ${\frak sp}_{8}$ & 
$\V_{\vp_4}$ & ${\frak f}_{4}$ & $\V_{\vp_1}$ & 3 \\ 
${\frak sl}_{2n+1}$ & --- & 
--- & ${\frak so}_{2n+1}$ & $\V_{2\vp_1}$ & --- \\ 
\hline
\end{tabular}
}  
\end{table}
\noindent
It follows from \re{look} that $\Delta_{\ov{0}}\cup\Delta_{\ov{1}}=
\Delta_{\ov{0}}$. Thus, everything lies in $\Delta_{\ov{0}}$.
Therefore a choice of the set of 
positive roots $\Delta_{\ov{0}}^+$
determines $\Delta_{\ov{1}}^+$, $\Delta_0^+$, and $\Delta_1^+$ as well. 
Of course, we choose $\Delta_{\ov{0}}^+$ so that
it is the image of $\Delta^+$ under the projection $\te^*\rar(\te_0)^*$.
\\
Then $\{\Delta_0^+,\Delta_1^+\}$ and $\{\Delta_{\ov{0}}^+,
\Delta_{\ov{1}}^+\}$ are two presentations for the totality of
all restricted {\sl positive\/} roots of $\g$.
Set $\rho_i=\frac{1}{2}|\Delta^+_i|$ and 
$\rho_{\ov{\imath}}=\frac{1}{2}|\Delta^+_{\ov{\imath}}|$
($i=0,1$).
Recall that $\rho=\frac{1}{2}|\Delta^+|$ and therefore $\Theta\rho=\rho$.
That is, $\rho\in (\te^*)^\Theta\simeq \te_0^*$.
It then follows from the above discussion that 
\beq \label{rho}
\rho=\rho_0+\rho_1=\rho_{\ov{0}}+\rho_{\ov{1}}=\rho_{\ov{0}}+
(\rho_{\ov{0}})_s\ .
\eeq
Let $W_0$ and $W_{\ov{0}}$ be the Weyl groups of $\g_0$ and $\g_{\ov{0}}$,
respectively. Although $\Delta_0\subset\Delta_{\ov{0}}$, $\g_0$ is not a
subalgebra of $\g_{\ov{0}}$ (if $\Theta\ne\OT$). In other words, $\Delta_0$
is a non-closed subset of $\Delta_{\ov{0}}$. Nevertheless, Prop.~\ref{coset}
applies to $W_0\subset W_{\ov{0}}$ and one obtains the subset $W'\subset
W_{\ov{0}}$ consisting of the elements of minimal length in the cosets
$\{wW_0\}$. Equivalently,
$W'=\{w\in W_{\ov{0}}\mid w(\Delta^+_0)\subset \Delta^+_{\ov{0}}\}$.
Below, we consider the Weyl chambers $C_0$ and $C_{\ov{0}}$, and
the hyperplanes $H_\mu$ ($\mu\in\Delta_1$).
They are regarded as subsets of the rational span of $\Delta_{\ov{0}}$
in $\te_0^*$.
\begin{s}{Proposition} \label{extr-outer} \\
1. The set of hyperplanes $H_\mu$ ($\mu\in\Delta_1$) cuts $C_0$ in 
$\# (W_{\ov{0}}/W_0)$ small chambers;  \\
2. The collection of ($\g_0$-dominant) weights 
$w^{-1}\rho_{\ov{0}}$ ($w\in W'$)
contains representatives of all small chambers in $C_0$. \\
3. The extreme weight of $Spin_0(\g_1)$ corresponding to 
$w^{-1}\rho_{\ov{0}}$ is
$\lb_w:=w^{-1}(\rho_0+\rho_1)-\rho_0=w^{-1}\rho-\rho_0$.
\end{s}\begin{proof}
To a great extent, the  proof is parallel to the proof of 
Prop.~\ref{inn-extr}.\\ 
1. The union $\Delta_0\cup\Delta_1$ coincides with $\Delta_{\ov{0}}$.
Therefore each small chamber is isomorphic to $C_{\ov{0}}$. 
Comparing the total number of chambers, we see that $C_0$ splits into 
$\# (W_{\ov{0}}/W_0)$ small chambers. \par
2 \& 3. By the definition of $W'$, it follows that $w^{-1}\rho_{\ov{0}}$
is $\g_0$-dominant. So, we have the required number of dominant weights
and it suffices to verify that the corresponding extreme weights of 
$Spin_0(\g_1)$ are different. 
\par
Given $w\in W'$,
the dominant half of $\Delta_1$ associated with $w^{-1}\rho_{\ov{0}}$ is 
\[
(\Delta_1)_{w}^+:=\{
\mu\in\Delta_1\mid (w^{-1}\rho_{\ov{0}},\mu)>0\}=\{
\mu\in\Delta_1\mid w\mu\in\Delta^+_{\ov{0}} \} 
\] 
and the corresponding extreme weight is
$\lb_w:=\frac{1}{2}|(\Delta_1)_{w}^+|$. 
\\ Set
$M_w=\{\mu\in\Delta_1^+\mid w\mu\in\Delta^+_{\ov{0}}\}$ and 
$\ov{M}_w=\Delta_1^+\setminus M_w$. Then 
\[
\rho_{\ov{0}}+\rho_{\ov{1}}=\rho_0+\rho_1=\rho_0+ \textstyle
\frac{1}{2}|M_w|+\frac{1}{2}|\ov{M}_w| \ .
\]
Since $w\in W'$, we have
\[
w^{-1}(\rho_{\ov{0}}+\rho_{\ov{1}})=\rho_0+ \textstyle
\frac{1}{2}|M_w|-\frac{1}{2}|\ov{M}_w| \ .
\]
Noting that $|(\Delta_1)_{w}^+|=|M_w|-|\ov{M}_w|$,
we obtain
$\lb_{w}=w^{-1}(\rho_{\ov{0}}+\rho_{\ov{1}})-\rho_0$. 
Obviously, these weights 
are different, and we are done.
\end{proof}%
In the next theorem, $\V_\lb$ denotes a $\g_0$-module.

\begin{s}{Theorem} \label{main-outer}
Let $\g=\g_0\oplus\g_1$ be a ${\Bbb Z}_2$-grading of outer type
and $\g=\g_{\ov{0}}\oplus\g_{\ov{1}}$ the associated diagram
${\Bbb Z}_2$-grading. Let
$W'$ be the set of representatives of minimal length for
$W_{\ov{0}}/W_0$. Then  \\[.6ex]
\centerline{
$Spin_0(\g_1)=\bigoplus_{w\in W'}\V_{\lb_w}$ .
} \vskip-.5ex
\end{s}\begin{proof}
First, note that if $\Theta$ is a diagram involutory automorphism,
then $W_{\ov{0}}=W_0$. Here the theorem claims that $Spin_0(\g_1)$ is 
irreducible, with highest weight $(\rho_{\ov{0}})_s$.
This was already demonstrated in Theorem~\ref{short}
and Prop.~\ref{B_n}. In the general case, we proceed as follows.
\\[.6ex]
By Prop.~\ref{extr-outer}, 
$\bigoplus_{w\in W'}\V_{\lb_w}\subset Spin_0(\g_1)$, and the equality
will follow from the fact that
$\dim (Spin_0(\g_1)^{\otimes 2})^{\g_0}=\# W'$. 
As in the proof of Theorem~\ref{main-inner}, a crucial step in the next
argument is of ``cohomological'' nature. Since
$\te_0$ contains regular elements, $\dim (\g_1)^{\te_0}=\dim\te-\dim\te_0$,
i.e., the multiplicity of the zero weight in $\g_1$
is equal to $\rk\g-\rk\g_0$. 
By Prop.~\ref{easy}(ii),
\[
\wedge^\bullet \g_1\simeq 
2^{\mbox{\scriptsize rk}\g-\mbox{\scriptsize rk}\g_0}{\cdot} 
Spin_0(\g_1)^{\otimes 2}
\]
and hence
\[
\dim(\wedge^\bullet \g_1)^{\g_0}= 
2^{\mbox{\scriptsize rk}\g-\mbox{\scriptsize rk}\g_0}{\cdot} 
\dim (Spin_0(\g_1)^{\otimes 2})^{\g_0} \ .
\]
At the rest of the proof, $\Bbbk=\Bbb C$.
Inspecting the list of the symmetric spaces of outer type and their 
cohomology rings 
over $\Bbb C$ (see e.g. \cite[\S\,4]{japan}) yields the equality
\[
\dim H^*(G/G_0)=
2^{\mbox{\scriptsize rk}\g-\mbox{\scriptsize rk}\g_0}{\cdot}
\# (W_{\ov{0}}/W_0)\ . 
\]
Since $(\wedge^\bullet\g_1)^{\g_0}\simeq H^*(G/G_0)$, 
we are done.
\end{proof}%
The following proof, although also being not free of case-by-case
arguments, does not appeal to $\Bbb C$. 

\begin{subs}{Another proof of Theorem~\ref{main-outer}} \label{another-outer}
\end{subs}
Arguing as in \re{another-inn} and using Prop.~\ref{extr-outer}(3), we obtain
\beq \label{tojdestvo}
\ch\Bigl(\bigoplus_{w\in W'}\V_{\lb_w}\Bigr)=\sum_{w\in W'}
\frac{\sum_{\tilde w\in W_0}\esi_0(\tilde w)
e^{\tilde w(\rho_0+\lb_w)}}{
\prod_{\ap\in\Delta_0^+}(e^{\ap/2}-e^{-\ap/2}) }= 
\eeq
\[
\phantom{\hspace{2cm}}=\frac{\sum_{w\in W'}
\sum_{\tilde w\in W_0}\esi_0(\tilde w)e^{\tilde w{w^{-1}}\rho}}{
\prod_{\ap\in\Delta_0^+}(e^{\ap/2}-e^{-\ap/2}) }= 
\frac{\sum_{w\in W_{\ov{0}}}\tau(w)e^{w\rho}}{
\prod_{\ap\in\Delta_0^+}(e^{\ap/2}-e^{-\ap/2}) } .
\]
Here $\tau(w)$ is the cunning parity for $w\in W_{\ov{0}}$,
relative to the subgroup $W_0$. To get another expression for the 
numerator, we exploit the following observation concerning the pairs
$(\g_0,\g_{\ov{0}})$ in Table~\ref{thety}.
Although $\Delta_0$ is not closed in $\Delta_{\ov{0}}$, the dual root
system $\widetilde\Delta_0$ is closed in $\widetilde\Delta_{\ov{0}}$ and,
moreover, it is a ``symmetric'' subset. That is, 
the partition $\widetilde\Delta_{\ov{0}}=\widetilde\Delta_{0}\sqcup
(\widetilde\Delta_{\ov{0}}\setminus \widetilde\Delta_{0})$ arises from an
{\it inner\/}
involutory automorphism of the ``dual'' Lie algebra. (E.g. the pair
$(\GR{C}{4},\GR{F}{4})$ inverts in $(\GR{B}{4},\GR{F}{4}$).)
Here $\widetilde\Delta_{\ov{0}}=(\Delta_{\ov{0}})_l \sqcup
2(\Delta_{\ov{0}})_s$. Recall from \re{look} that
$\Delta_{\ov{1}}=(\Delta_{\ov{0}})_s=(\Delta_0)_s$. 
This means in particular that
$\Delta_{\ov{0}}\setminus \Delta_0$ consists of long roots and these 
are exactly the roots constituting $(\Delta_1)_l$. 
Hence $\widetilde\Delta_{\ov{0}}\setminus \widetilde\Delta_{0}
=(\Delta_1)_l$. After these preparations, write out the identity
from Theorem~\ref{identity} for the partition
$\widetilde\Delta_{\ov{0}}=\widetilde\Delta_{0}\sqcup (\Delta_1)_l$:
\beq \label{dopolnit}
\sum_{w\in W_{\ov{0}}}\tau(w)e^{w\tilde\rho}=
\prod_{\ap\in\widetilde\Delta_0^+}(e^{\ap/2}-e^{-\ap/2})
\prod_{\mu\in (\Delta_1^+)_l}(e^{\mu/2}+e^{-\mu/2}) .
\eeq
Here $\tilde\rho:=\frac{1}{2}|\widetilde\Delta_{\ov{0}}^+|=
\frac{1}{2}|(\Delta_{\ov{0}}^+)_l| + |(\Delta_{\ov{0}}^+)_s|=
(\rho_{\ov{0}})_l+2(\rho_{\ov{0}})_s=\rho_{\ov{0}}+(\rho_{\ov{0}})_s=
\rho$ (see Eq.~\ref{rho}). Transforming the first factor
on the right hand side of Eq.~\re{dopolnit} yields
\[
\prod_{\ap\in\widetilde\Delta_0^+}(e^{\ap/2}-e^{-\ap/2})=
\prod_{\ap\in(\Delta_0^+)_s}(e^{\ap}-e^{-\ap})
\prod_{\beta\in(\Delta_0^+)_l}(e^{\beta/2}-e^{-\beta/2})=
\]
\[
=\prod_{\ap\in(\Delta_0^+)_s}(e^{\ap/2}-e^{-\ap/2})(e^{\ap/2}+e^{-\ap/2})
\prod_{\beta\in(\Delta_0^+)_l}(e^{\beta/2}-e^{-\beta/2})=
\]
\[
=\prod_{\ap\in\Delta_0^+}(e^{\ap/2}-e^{-\ap/2})
\prod_{\mu\in(\Delta_0^+)_s}(e^{\mu/2}+e^{-\mu/2}) .\phantom{\hspace{3cm}}
\]
Since $(\Delta_0^+)_s=(\Delta_1^+)_s$, the whole expression on
the right hand side of Eq.~\re{dopolnit}  is equal to
\[
\prod_{\ap\in\Delta_0^+}(e^{\ap/2}-e^{-\ap/2})
\prod_{\ap\in\Delta_1^+}(e^{\mu/2}+e^{-\mu/2}) .
\]
Substituting this expression for 
$\sum_{w\in W_{\ov{0}}}\tau(w)e^{w\rho}$ in Eq.~\re{tojdestvo},
we obtain
\[ \phantom{\hspace{3cm}}
\ch\Bigl(\bigoplus_{w\in W'}\V_{\lb_w}\Bigr)=
\prod_{\ap\in\Delta_1^+}(e^{\mu/2}+e^{-\mu/2})=\ch Spin(\g_1) .
\phantom{\hspace{2cm}} \mbox{q.e.d.}
\]
\begin{subs}{Connexion with cohomology of symmetric spaces}
\label{deeper}\end{subs}
Using classical structure results on $H^*(G/G_0)$ 
(see e.g. \cite[ch.3]{al}), one may notice some 
interesting coincidences. \\
{\sf 1}. Recall that the image of the canonical map $\eta : H^*(G/G_0)\rar
H^*(G)$ is generated by primitive elements.
Actually, there exists a subalgebra of
$H^*(G/G_0)$ that is mapped isomorphically onto $\Ima(\eta)$.
It is called {\it Samelson\/} and denoted by $Sam(G/G_0)$.
There exists also another subalgebra of $H^*(G/G_0)$, which is called
{\it characteristic\/} and denoted by $_0H^*(G/G_0)$. It is the zero component
for a natural grading in $H^*(G/G_0)$. Then
$H^*(G/G_0)\simeq Sam(G/G_0)\otimes\, _0H^*(G/G_0)$. This holds not only
for the symmetric spaces but for a wider class of {\it formal\/}
homogeneous spaces, see [loc.\,cit., \S\,12,\,Th.2].
A general fact for the formal homogeneous spaces
is that $\dim Sam(G/G_0)=
2^{\mbox{\scriptsize rk}\g-\mbox{\scriptsize rk}\g_0}$. Hence
$\dim Spin_0(\g_1)=\dim\,_0H^*(G/G_0)$, which suggests that
$Spin_0(\g_1)$ might somehow be related to the characteristic subalgebra of
$H^*(G/G_0)$. \\
{\sf 2}. Note that $H^*(G/G_0)=\, _0H^*(G/G_0)$ if and only if 
$\Theta$ is inner,
and  $H^*(G/G_0)=Sam(G/G_0)$ if and only if $\Theta$ is
a diagram involutory automorphism. In the mixed case, the associated
diagram involutory automorphism $\OT$ seems to yield a splitting for
$H^*(G/G_0)$. Namely, one has
$\dim Sam(G/G_0)=\dim H^*(G/G_{\ov{0}})$ and
$\dim\, _0H^*(G/G_0)=\dim H^*(\widetilde G_{\ov{0}}/\widetilde G_0)$, where
$\widetilde G_{\ov{0}}$ and $\widetilde G_0$ are the groups corresponding to 
the dual root systems $\widetilde\Delta_{\ov{0}}$ and
$\widetilde\Delta_{{0}}$. Moreover, these two pairs of {\sl graded\/} algebras
have equal Poincar\'e polynomials and it is likely that they are naturally
isomorphic. I think that a better understanding of this situation as
well as elimination of the
case-by-case arguments in section~\ref{spin5} can be achieved 
through application of the theory of twisted affine Kac--Moody algebras.
\\
{\sf 3}. In the above exposition, $\OT$ has appeared as {\sl deus ex machina}.
But the Kac-Moody theory provides some explanation for this. Namely,
the outer automorphism $\Theta$ determines the twisted affine Kac-Moody
algebra $\hat{{\cal L}}(\g,\Theta,2)=\hat{{\cal L}}(\g)$ \cite[ch.\,8]{kac}.
Next, $\hat{{\cal L}}(\g)$ has the standard ${\Bbb Z}$-grading 
associated with the special vertex of the Dynkin diagram of 
$\hat{{\cal L}}(\g)$. If $\g\ne {\frak so}_{2n+1}$, this $\Bbb Z$-grading 
determines
another outer automorphism of $\g$, which is just $\OT$. 
\\[.7ex]
\refstepcounter{equation}\label{ex-outer}%
{\bf \theequation\quad Examples.}\quad
{\sf 1}. $\g={\frak sl}_{2n}$, $\g_0={\frak so}_{2n}$.
Then $\g_1\simeq\V_{2\vp_1}$. As indicated in Table~\ref{thety},
$\g_{\ov{0}}={\frak sp}_{2n}$ and therefore $\# (W_{\ov{0}}/W_0)=2$. 
Here $\Delta_0=\{\pm\esi_i\pm\esi_j\mid 1\le i,j\le n, i\ne j\}$
and $\Delta_{\ov{0}}=\{\pm\esi_i\pm\esi_j, \pm 2\esi_i\}$. This example is
a kind of outer version of Example~\ref{ex-inner}(1). Indeed,
taking the ``dual'' Lie
algebras for $(\g_0,\g_{\ov{0}})$ yields the symmetric pair
considered there.
In our case, $\Delta_1=\Delta_{\ov{0}}$ and
$W'=\{id, w_n\}$, where $w_n(\esi_i)=\esi_i$ ($i\le n-1$)
and $w_n(\esi_n)=-\esi_n$. Therefore
$(\Delta_1)^+_{id}=\Delta_{\ov{0}}^+=\Delta_0^+\cup\{2\esi_1,\dots,2\esi_n\}$
and
$(\Delta_1)^+_{w_n}=\Delta_0^+\cup\{2\esi_1,\ldots,2\esi_{n-1},-2\esi_n\}$. 
Hence the extreme weights are $\rho_0+2\vp_n$ and $\rho_0+2\vp_{n-1}$.
Thus, $Spin_0(\V_{2\vp_1})=
\V_{\rho_0+2\vp_{n-1}}\oplus\V_{\rho+2\vp_n}$ and, 
because $m_{2\vp_1}(0)=n-1$,
$\wedge^\bullet\V_{2\vp_1}=2^{n-1}{\cdot}
(\V_{\rho+2\vp_{n-1}}\oplus\V_{\rho+2\vp_n})^{\otimes 2}$.
\\[.6ex]
{\sf 2}. $\g={\frak e}_6$, $\g_0={\frak sp}_8$. 
Then $\g_1\simeq\V_{\vp_4}$. As indicated in Table~\ref{thety}, 
$\g_{\ov{0}}={\frak f}_4$ and hence
$\# (W_{\ov{0}}/W_0)=3$. 
This is the outer version of Example~\ref{ex-inner}(2).
Here
$\Delta_{\ov{0}}=\{\pm\esi_i\pm\esi_j, \pm\esi_i, (\pm\esi_1\pm\esi_2
\pm\esi_3\pm\esi_4)/2\mid 1\le i,j\le 4, i\ne j\}$ and
$\Delta_0=(\Delta_{\ov{0}})_s\cup \{\pm\esi_1\pm\esi_2,\,\pm\esi_3\pm\esi_4
\}$. The standard set of simple roots for $\Delta_{\ov{0}}$ is
$\ap'_1=\frac{1}{2}(\esi_1-\esi_2-\esi_3-\esi_4)$,
$\ap'_2=\esi_4$,
$\ap'_3=\esi_3-\esi_4$,
$\ap'_4=\esi_2-\esi_4$.
The roots for ${\frak sp}_8$ have non standard presentation, but it is not
hard to find that the simple roots in $\Delta^+_{\ov{0}}\cap\Delta_0$ are
$\ap_1=\esi_2$,
$\ap_2=\frac{1}{2}(\esi_1-\esi_2-\esi_3-\esi_4)$,
$\ap_3=\esi_4$,
$\ap_4=\esi_3-\esi_4$. Therefore the fundamental weights for ${\frak sp}_8$
are
$\vp_1=\frac{1}{2}(\esi_1+\esi_2)$,
$\vp_2=\esi_1$,
$\vp_3=\esi_1+\frac{1}{2}(\esi_3+\esi_4)$,
$\vp_4=\esi_1+\esi_3$.\\
Then an explicit verification shows that $W'=\{id,w',w''\}$, where
$w',w''$ are permutations of $\{\esi_1,\esi_2,\esi_3,\esi_4\}$ determined by
the cycles (23) and (432), respectively. (E.g. $w''(\esi_4)=\esi_3$.)
The direct computation of the extreme weights gives:

$\lb_{id}=\rho_{\ov{0}}+\rho_{\ov{1}}-\rho_0=
\frac{1}{2}(9\esi_1+5\esi_2+\esi_3+\esi_4)=5\vp_1+\vp_2+\vp_3$,

$\lb_{w'}=(w')^{-1}(\rho_{\ov{0}}+\rho_{\ov{1}})-\rho_0=
\frac{1}{2}(9\esi_1+3\esi_2+3\esi_3+\esi_4)=3\vp_1+\vp_2+\vp_3+\vp_4=
\rho_0+2\vp_1$,

$\lb_{w''}=(w'')^{-1}(\rho_{\ov{0}}+\rho_{\ov{1}})-\rho_0=
\frac{1}{2}(9\esi_1+\esi_2+3\esi_3+3\esi_4)=\vp_1+\vp_2+3\vp_3$.
\\
Whence \\
\centerline{ $Spin_0(\V_{\vp_4})=\V_{5\vp_1+\vp_2+\vp_3}
\oplus\V_{\vp_1+\vp_2+3\vp_3}\oplus
\V_{\rho_0+2\vp_1}$}
\\ and
\[
\wedge^\bullet\V_{\vp_4}=4\,(\V_{5\vp_1+\vp_2+\vp_3}
\oplus\V_{\vp_1+\vp_2+3\vp_3}\oplus
\V_{\rho_0+2\vp_1})^{\otimes 2} \ .
\]

\sekt{Decomposably-generated 
`Spin' modules and  a Casimir element\label{spin6}}%
Recall that we have given a geometric description of 
the extreme weights of Spin-representations in \re{extreme}.
\\
{\bf Definition.}  Given an orthogonal $\g$-module $\V$, 
the $\g$-submodule of $Spin_0(\V)$ generated by the extreme weight vectors
is denoted by $Spin^{dg}_0(\V)$; $Spin_0(\V)$ 
is called {\it decomposably-generated\/}, if it is equal to 
$Spin^{dg}_0(\V)$, i.e., if all its highest weights are extreme. 
\\[1ex]
Since the extreme weights are of multiplicity 1,
``decomposably-generated'' implies ``multiplicity free''.
As a consequence of previous development, we have
\begin{s}{Proposition} \label{dg}
Let $\g=\g_0\oplus\g_1$ be a ${\Bbb Z}_2$-graded semisimple Lie algebra. Then
the $\g_0$-module $Spin_0(\g_1)$ is decomposably-generated (and
multiplicity free). 
\end{s}\begin{proof}
The problem immediately reduces to the case in which $\g$ is an irreducible
${\Bbb Z}_2$-graded algebra. Then either $\g$ is simple or 
$\g\simeq\h\oplus\h$, where $\h$ is simple and $\Theta(h_1,h_2)=(h_2,h_1)$.
In the second case, $\g_0\simeq\h$ is the diagonal in $\g$, and $\g_1\simeq
\g_0$ as $\g_0$-module. Here the conclusion follows by Kostant's result,
see Example~\ref{ex-adj}(1). In the first case,
for $\g_0$ of inner type, use
Prop.~\ref{inn-extr}(3) and Theorem~\ref{main-inner};
for $\g_0$ of outer type, use
Prop.~\ref{extr-outer}(3) and Theorem~\ref{main-outer}.
\end{proof}%
An explanation of the term ``decomposably-generated'' comes from 
Example~\ref{ex-adj}(3). If $\V=\W\oplus\W^*$, then $Spin(\V)\simeq
\Bbbk_{-\nu}\otimes\wedge^\bullet\W$ and each extreme weight vector is 
represented by a decomposable vector in the exterior algebra.
\\
I think that the property of being ``decomposably-generated'' 
characterizes the representations of the form $Spin_0(\g_1)$, i.e.,
\begin{s}{Conjecture} \label{gipoteza}
Let $\V$ be an orthogonal $\g$-module. Then
$Spin_0(\V)$ is decompo\-s\-ab\-ly-generated if and only if 
$\tilde\g:=\g\oplus\V$ is
a ${\Bbb Z}_2$-graded semisimple Lie algebra.
\end{s}%
The conjecture will be proved in a particular case.
Until the end of the section, the following situation is being considered:
$\g$ is semisimple, $\h$ is a reductive subalgebra of $\g$, 
and $\me:=\h^\perp\subset\g$.
Then $\me$ is an orthogonal $\h$-module and $\g=\h\oplus\me$ is a vector
space sum. Clearly, this decomposition is a ${\Bbb Z}_2$-grading if and only
if $[\me,\me]\subset\h$.
The representation $\h\rar{\frak so}(\me)$ is the isotropy representation of
the affine homogeneous space $G/H$. Our aim is to study
$Spin_0(\me)$ and $Spin_0^{dg}(\me)$ in the equal rank  case. That is,
it is assumed from now on that $\rk\h=\rk\g$ and, more precisely,
$\te\subset\h$.
Then $\Delta^+=\Delta_\h^+\sqcup\Delta_\me^+$, $\me$ has no zero weight and
$\wedge^\bullet\me\simeq (Spin_0(\me))^{\otimes 2}$.
Denoting by $W_\h$ the Weyl group of $(\h,\te)$, one may consider the 
minimal length ``section'' $W^\h$ for $W\rar W/W_\h$ and the cunning
parity $\tau: W\rar\{1,-1\}$, determined by $\Delta_\h^+$. 
The proof of Prop.~\ref{inn-extr} applies in the present situation as well.
This yields exactly $\# W^\h$ extreme weights of $Spin_0(\me)$. Hence
\[
\dim (\wedge^\bullet\me)^\h=\dim (Spin_0(\me)^{\otimes 2})^\h\ge \# W^\h \ .
\]
For $w\in W^\h$, the corresponding extreme weight is $\lb_w=
w^{-1}\rho-\rho_\h$, where $\rho_\h=\frac{1}{2}|\Delta^+_\h|$. Therefore,
arguing as in \re{another-inn}, we obtain
\[
\ch(Spin^{dg}_0(\me))=
\ch\Bigl(\bigoplus_{w\in W^\h}\V_{\lb_w}\Bigr) =
\frac{\sum_{w\in W}\tau(w)e^{w\rho}}{
\prod_{\ap\in\Delta^+_\h}(e^{\ap/2}-e^{-\ap/2})} \ .
\]
Recall that $\ch Spin_0(\me)=
\prod_{\mu\in\Delta^+_\me}(e^{\mu/2}+e^{-\mu/2})$.
On the other hand, $\dim H^*(G/H)=\# W^\h$ 
\cite[\S\,13,\,Th.\,2]{al} and $H^*(G/H)$ can be computed via the
complex of $G$-invariant exterior forms on $G/H$, i.e., the complex
$\bigl( (\wedge^\bullet(\g/\h)^*)^\h, d\bigr)$, where $d$ is the usual
Lie algebra coboundary operator.  We shall identify
the $\h$-modules $(\g/\h)^*$ and $\me$.
Having compared the previous expressions, we obtain
\begin{s}{Proposition}  \label{4-usloviya}
Let $\h\subset\g$ be a reductive subalgebra of maximal rank and $\h\rar{\frak so}(\me)$
the isotropy representation. Then
the following conditions are equivalent:

{\sf (i)} \ $d$ is trivial on $(\wedge^\bullet\me)^\h$;

{\sf (ii)} $\dim (\wedge^\bullet\me)^\h= \# W^\h$;

{\sf (iii)} $Spin_0(\me)=Spin^{dg}_0(\me)$;

{\sf (iv)} $\sum_{w\in W}\tau(w)e^{w\rho}=
\prod_{\ap\in\Delta^+_\h}(e^{\ap/2}-e^{-\ap/2})
\prod_{\mu\in\Delta^+_\me}(e^{\mu/2}+e^{-\mu/2})$. \qu
\end{s}%
Thus, conjecture~\ref{gipoteza}
claims that neither of these conditions holds
unless $G/H$ is symmetric.
\begin{s}{Theorem}  \label{d3}
Let $\h$ be a reductive subalgebra of $\g$, with $\rk\h=\rk\g$, and
$\me:=\h^\perp$. Suppose $[\me,\me]\not\subset\h$; then  $d$
is non-trivial on $(\wedge^\bullet\me)^\h$. More precisely,
$d((\wedge^3\me)^\h)\ne 0$.
\end{s}\begin{proof*} 
For any $x\in\g$, let $x_\h$ and $x_\me$ denote its components in $\h$ 
and $\me$, respectively. Given $x,y\in\me$, consider the decomposition
$[x,y]=[x,y]_\h+[x,y]_\me$. We regard $[\ ,\ ]_\me$ as mapping from
$\me\times\me$ to $\me$, and likewise for $[\ ,\ ]_\h$.
By assumption, $[\ ,\ ]_\me\not\equiv 0$. On the other hand, applying
construction from section~\ref{grassmann} (see~\re{mubar} and around)
to $\h$ and $\me$ in place of $\g$ and $\V$, we see that
$[x,y]_\h=\bar\mu(x,y)$ for any $x,y\in\me$.

Define the 3-form $\Psi: \wedge^3\me\rar\Bbbk$ by
$\Psi(x,y,z)=\Phi([x,y],z)$. Obviously, $\Psi$ is $\h$-invariant.
The assumption $[\me,\me]\not\subset\h$ precisely means that $\Psi
\not\equiv 0$. We shall prove $d\Psi\ne 0$. To compute $d\Psi$, we regard
$\Psi$ as $\h$-invariant 3-form on $\g$, orthogonal to $\h$,
and use the standard formula for $d$. The resulting expression 
is 
\[
d\Psi(x,y,z,u)=
2\Bigl( \Phi([x,y]_\me,[z,u]_\me)+\Phi([y,z]_\me,[x,u]_\me)+
\Phi([z,x]_\me,[y,u]_\me)\Bigr) \ .
\]
Since $[x,y]=\bar\mu(x,y)+[x,y]_\me$, $\h$ is orthogonal to $\me$, and\\
$([x,y],[z,u])+([y,z],[x,u])+([z,x],[y,u])=0$ 
(because of the Jacobi identity), we have
\[
d\Psi(x,y,z,u)=
-2\Bigl( \Phi(\bar\mu(x,y),\bar\mu(z,u))+\Phi(\bar\mu(y,z),\bar\mu(x,u))+
\Phi(\bar\mu(z,x),\bar\mu(y,u))\Bigr) \ .
\]
In the notation of Prop.~\ref{jacobi}, for $\h$ and $\me$ in place of
$\g$ and $\V$, this means that $d\Psi=-2\kappa$.
Assume that $\kappa=0$. Then $\h\oplus\me$ equipped with the modified
multiplication $[\ ,\ ]\widetilde{\ }$ becomes a ${\Bbb Z}_2$-graded
Lie algebra (see Prop.~\ref{jacobi}). Clearly, the multiplication
does change only for pairs of elements in $\me$: $[m_1,m_2]\widetilde{\ }:
=[m_1,m_2]_\h$, whereas the structure of Lie algebra on $\h$ and the
$\h$-module structure on $\me$ remain undisturbed.
Let $\tilde\g$ denote the Lie algebra with modified multiplication
and $\Theta$ the corresponding involutory automorphism of $\tilde\g$.
It is easily seen that $\tilde\g$ is semisimple (use the proof of
Theorem~\ref{main1}) and $\Theta$ is inner (because $\te\subset\h$
remains a Cartan subalgebra in $\tilde\g$).
For the symmetric space $\tilde G/H$, we have
$H^3(\tilde G/H)=(\wedge^3\me)^\h\ne 0$. But $H^{odd}({\cdot})=0$
for the symmetric spaces of inner type \cite[\S\,13,\,n.3]{al}. 
This contradiction proves $\kappa=d\Psi\ne 0$.
\end{proof*}%
\begin{s}{Corollary}
Conjecture \ref{gipoteza} is true for the isotropy representations
of affine homogeneous spaces $G/H$ with $\rk\g=\rk\h$. \qu
\end{s}%
{\bf Remark.} Theorem~\ref{d3} is true even if $\rk\h<\rk\g$ and some mild
conditions are satisfied (e.g. $\g$ is simple). However this has no immediate
relation to Conjecture~\ref{gipoteza}.
\\[.8ex]
For a reductive Lie algebra $\h$, the Casimir element in $U(\h)$
is determined by the choice of an invariant bilinear form on $\h$.
If $\h$ is not simple, then the choice is essentially non
unique. But for the isotropy representations
one has a preferred choice of the bilinear form. In the above setting,
let $\Phi(\ ,\ )_\h$ be the restriction of $\Phi(\ ,\ )$ to $\h$.
Notice that even if $\h$ is semisimple and we begin with the Killing form
on $\g$, then$\Phi(\ ,\ )_\h$ is not necessarily
proportional to the Killing form on $\h$.
Let $c_\h\in U(\h)$ be the Casimir element with respect to $\Phi(\ ,\ )_\h$. 
Recall that the $W$-invariant
scalar product on ${\cal P}_{\Bbb Q}$ is determined by $\Phi(\ ,\ )$.
\begin{s}{Proposition}  \label{casimir}
Suppose $\rk\h=\rk\g$. Then
the Casimir element $c_\h$ acts scalarly on $Spin^{dg}_0(\me)$.
Its eigenvalue is equal to $(\rho,\rho)-(\rho_\h,\rho_\h)$. 
\end{s}\begin{proof*}
As is indicated above, $Spin_0^{dg}(\me)=\oplus_{w\in W^\h}\V_{\lb_w}$
and $\lb_w=w^{-1}\rho-\rho_\h$. Therefore the eigenvalue of $c_\h$ on 
$\V_{\lb_w}$ is $(\lb_w+2\rho_\h,\lb_w)=(w^{-1}\rho,w^{-1}\rho)-
(\rho_\h,\rho_h)=(\rho,\rho)-(\rho_\h,\rho_\h)$. 
\end{proof*}%
\begin{s}{Theorem}  \label{inv-casimir}
Let $\g=\g_0\oplus\g_1$ be a ${\Bbb Z}_2$-graded semisimple Lie algebra.
Define the Casimir element $c_0$ for $\g_0$ using the restriction of
$\Phi(\ ,\ )$ to $\g_0$. Then $c_0$ acts on
$Spin_0(\g_1)$ scalarly, with value $(\rho,\rho)-(\rho_0,\rho_0)$.
\end{s}%
\begin{proof}
1. If $\Theta$ is inner, then $\rk\g_0=\rk\g$; we conclude by Propositions
\ref{dg}, \ref{casimir}.\\
2. If $\Theta$ is outer, some accuracy is needed, since
$\rk\g_0< \rk\g$. 
We use notation and information from section~\ref{spin5}.
Since $Spin_0(\g_1)=Spin^{dg}_0(\g_1)=\oplus_{w\in W'}\V_{\lb_w}$ and
$\lb_w=w^{-1}\rho-\rho_0$, the value of $c_0$ on
$\V_{\lb_w}$ is equal to 
$(w^{-1}\rho,w^{-1}\rho)-(\rho_0,\rho_0)$. Here $W'\subset W_{\ov{0}}$, 
where $W_{\ov{0}}$ is the Weyl group of $\g_{\ov{0}}$. Recall that 
$\te_0=\te^\Theta$ is a Cartan subalgebra for both
$\g_{\ov{0}}$ and $\g_0$.  
As $\te_0$ contains regular elements of $\te$, we have
$N_{G_{\ov{0}}}(\te_0)
\subset N_G(\te)$. Furthermore, since $G_{\ov{0}}$ is connected
and $\te_0$ is Cartan, we have $G_{\ov{0}}\cap T=T_0$.
Therefore 
$W_{\ov{0}}=N_{G_{\ov{0}}}(\te_0)/T_0$
can be identified with a subgroup of $W=N_G(\te)/T$. Hence 
$(w^{-1}\rho,w^{-1}\rho)=(\rho,\rho)$.
\end{proof}

\vno{3} \indent
{\footnotesize 
{\bf }\ \parbox[t]{190pt}{%
 {\it {\hspace*{1ex}} \\ ul. Akad. Anokhina \\
d.30, kor.1, kv.7 \\
Moscow 117602 \quad Russia} \\ E-mail: dmitri@panyushev.mccme.ru 
} \qquad \parbox[t]{190pt}{%
Current address (until February 29, 2000): \\
{\it Max-Planck-Institut f\"ur Mathematik \\
Vivatsgasse 7 \\
D-53111 Bonn, Deutschland} \\ 
e-mail: panyush@mpim-bonn.mpg.de
}
}
\end{document}